\documentclass[11pt]{article}
\usepackage[letterpaper,margin=1in]{geometry}
\usepackage[T1]{fontenc}
\usepackage{lmodern}
\usepackage{microtype}
\usepackage[nocompress]{cite}
\usepackage{graphicx}
\usepackage{tikz}
\usetikzlibrary{fit}
\usepackage{amsmath, amssymb, amsfonts, amsthm}
\usepackage{algorithm}
\usepackage{algpseudocode}
\usepackage{mathtools}
\usepackage{multirow}
\usepackage{mathrsfs}
\usepackage{makecell}
\usepackage{bbding}
\usepackage{authblk}

\usepackage[hidelinks]{hyperref}
\usepackage{orcidlink}
\newtheorem{definition}{Definition}
\newtheorem{problem}{Problem}

\newtheorem{corollary}{Corollary}
\newtheorem{theorem}{Theorem}
\newtheorem{example}{Example}
\DeclarePairedDelimiter\floor{\lfloor}{\rfloor}

\title{A Column-Generation Approach to Length-Constrained Shortest Paths under Cascading Causal Failures}
\author[1,2]{Nicol\`o Aicardi\,\orcidlink{0009-0000-4796-4171}}
\author[3]{Zuyuan Zhang\,\orcidlink{0000-0001-5795-9771}}
\author[3]{Sridhar Radhakrishnan\,\orcidlink{0000-0002-0327-8134}}
\author[4]{Kash Barker\,\orcidlink{0000-0002-0142-1558}}
\author[1]{Dritan Nace\,\orcidlink{0000-0003-3914-2463}}
\affil[1]{Heudiasyc, UMR CNRS 7253, Universit\'e de Technologie de Compi\`egne, 57 avenue de Landshut, CS 60319, 60200 Compi\`egne, France}
\affil[2]{Department of Informatics, Robotics and Systems Engineering, Universit\`a degli Studi di Genova, Via Dodecaneso 35, 16146 Genova, Italy}
\affil[3]{Department of Computer Science, University of Oklahoma, Norman, OK 73019, USA (e-mail: zuyuan.zhang-1@ou.edu; sridhar@ou.edu)}
\affil[4]{School of Industrial and Systems Engineering, University of Oklahoma, Norman, OK 73019, USA (e-mail: kashbarker@ou.edu)}
\date{}

\begin{document}
\maketitle
{\renewcommand{\thefootnote}{}\footnotetext{This work has been submitted to the IEEE for possible publication. Copyright may be transferred without notice, after which this version may no longer be accessible.}}

\begin{abstract}
In the study of cascading causal node failures, the failure of one set of nodes can trigger the failure of another, possibly non-adjacent, set, and these triggers can propagate in cascades.
Given a graph, a set of causal failures, a source--terminal pair, and a constant $\alpha \geq 1$, we study the length-constrained shortest-path problem of finding the maximum subset of causalities whose application preserves a source--terminal path of length at most $\alpha$ times the original distance.
We show that the problem is NP-hard, a hardness that reflects the non-local, cascading nature of causal failures.
We propose a mixed-integer linear programming formulation with explicit node-integrity variables that decouple the node state from routing, a column-generation algorithm built on a path-based reformulation whose pricing subproblem is solved by a three-level cascade (beam search, label-capped $A^{*}$, and exact $A^{*}$ with warm start), and a greedy heuristic as a fast baseline.
On the IEEE~300 and Power~5000 benchmarks, all three methods are competitive.
In dense synthetic networks of $10{,}000$ to $15{,}000$ nodes, the exact model exhausts its time budget or returns an order of magnitude fewer causalities, whereas column generation remains fast and accurate, emerging as the method of choice for large instances.
\end{abstract}


\begin{table}[!t]
    \caption{Notation}
    \label{tab:notation}
    \centering
    \footnotesize
    \begin{tabular}{@{}l p{0.6\textwidth}@{}}
    \hline
        $G = (V, E)$ & An undirected graph  \\
        $V$          & The node set \\
        $V_i$        & The $i$-th subset of $V$\\
        $U$          & A subset of $V$ \\
        $E$          & The edge set \\
        $n$          & The number of nodes \\
        $m$          & The number of edges \\
        $C: V_{1}\Rightarrow V_2$ & A causality with trigger set $V_1$ and effect set $V_2$ (also written $V_{1,C}, V_{2,C}$; and $V_{1,i}, V_{2,i}$ for the $i$-th causality $C_i$)\\
        $\mathcal{C}$ & The set of all causalities\\
        $h=|\mathcal{C}|$          & The number of causalities\\
        $s, t$       & The source and terminal nodes\\
        $uv$-path    & A path from $u$ to $v$\\
        $\alpha$     & The length constraint factor of $s$--$t$ paths\\
        $D_0$        & The distance from $s$ to $t$ in $G$\\
        $D=\alpha\times D_0$         & The length constraint\\
        $d_G(u, v)$  & The distance between $u$ and $v$ in graph $G$\\
        $P$          & An $st-$path with length smaller than or equal to $D$\\
        $\mathcal{P}$ & The set of $st-$paths with length smaller than or equal to $D$, i.e. $P\in \mathcal{P}$\\
        MILP & The solution obtained by the optimization model\\
        CG & The solution obtained by the column generation framework\\
        Greedy & The solution obtained by heuristic algorithm\\
        \hline
    \end{tabular}
\end{table}
\section{Introduction}\label{sec:Introduction}
Network infrastructure systems, such as electricity, transportation, and communication networks, are of critical importance for the daily function of society. These interconnected systems underpin virtually every aspect of modern life, from powering homes and enabling commerce to facilitating the flow of information across the globe. As societies grow increasingly dependent on these networks, any disruption, whether caused by natural disasters, technical failures, or deliberate attacks~\cite{Rinaldi2001,Buldyrev2010}, can have far-reaching consequences for public safety, economic stability, and national security.

Infrastructure networks are inherently complex systems, characterized by intricate interactions and interdependencies among their constituent components (nodes/edges). Within these networks, individual components (e.g., power stations, transport hubs, communication relays) interact through the exchange of material, energy, or information, forming a tightly coupled web of dependencies. These interdependencies can give rise to emergent phenomena that introduce novel reliability challenges when failures occur. 

These challenges are especially impactful in the presence of causal failures~\cite{zzy_causal}, which imply that the failure of some nodes would cause the failure of some other nodes, and the failed nodes may not necessarily be adjacent to each other. As such, causal failures may be the result of non-obvious dependencies, further adding to the complexity of understanding their reliability. 

Our prior work explored the impact of causal failures in the form of the sizes of the remaining giant component and small components~\cite{zzy_causal,zzy_augment,nanyanzi_disinfo}. In this paper, we focus on a different dimension of network vulnerability: the impact of causal failures on the distance between a particular pair of nodes.  The following examples motivate the study of such causal failures and their impacts on the shortest paths in networks.

{\it Electrical Grid System}: We consider an electrical grid system with several generators (e.g., power stations, batteries), demand points as nodes, and transmission lines as edges connecting the nodes. 
Given that some generating stations share a load, the failure of one generating station may cause the failure of other generators.
For example, if the maximum load of a station, say node $i$, is beyond its capacity, the system may distribute the additional load to other stations, say node $j$~\cite{Motter2002,Dobson2007}. 
If node $j$ cannot meet the demand, it will fail. 
Therefore, we have the failure of node $j$ due to failure of node $i$.
For a given pair of nodes $s, t$, when causal failures occur, the distance from $s$ to $t$ might increase, so it takes longer to transmit power from $s$ to $t$ and generates more power loss \cite{Min}.

{\it Interdependent Infrastructure Systems} \cite{Sen}: Consider a system of interdependent power and communication networks with a smart control center within the system. 
The smart control center is operational only if it receives power from the generators of the power network. 
Suppose that there are only two disjoint paths between the generator $s$ and the control center $t$ (e.g., $\{s, a, b, t\}$ and $\{s, c, d, t\}$ where $a, b, c, d$ are intermediate nodes in the system).
Then one of $\{a, b\}$ and $\{c, d\}$ must be operational to make power available to the control center.
 In contrast, if nodes $a$ and $c$ fail, then the control center fails because no path is available.
Similarly, failures of the nodes $b$ and $d$ trigger the failure of the control center. 
Moreover, node failures in the system may increase the distance between a given pair of nodes in the communication network and thus harm the quality of communication.
For example, a fiber communication network uses the wavelength division multiplexing network, providing multiple channels over a single fiber employing different light wavelengths~\cite{Mukherjee2000}.
The long distance may cause network signal attenuation problems and node overload \cite{Mutafungwa:01}.
Hence, there is usually a path length constraint in communication networks.

{\it Networked Systems with Multiple Operating Systems}: We consider the situation where routers in one network run different versions of an operating system, where a cyberattack occurs on a specific version of the operating system on a router. 
This attack results in the exposure of other routers with that version of the operating system to attacks.
Hence, due to the cyberattack, the failure of one router triggers the failure of other routers.
Moreover, these routers work between a pair of workstations, say $s$ and $t$, to transmit data. 
Generally, if all routers function, the data would take the shortest possible time to finish transmission. 
However, once router failures occur, the number of hops for data transmission would potentially increase. 
More specifically, suppose that each router fails with the same probability $q=1-p$, then the probability that an $s$--$t$ path is operational is $p^l$, where $l$ is the number of routers~\cite{BarlowProschan1996}.
Once the number of hops increases, the probability is lower than before, resulting in the path being more likely to fail.

A common thread across these examples is that causal failures are considerably harder to analyze than localized node failures studied in classical reliability models due to the complexity they introduce. When a failure can propagate only to adjacent nodes, the consequences of a disruption can be reasoned about locally; causal failures remove this guarantee. 
Since the triggered nodes may be arbitrarily far from the nodes that fail and since a single failure may set off a cascade through the dependency structure, one must first compute the \emph{closure} of the failed set before its effect on connectivity and distances can be evaluated. This non-local, cascading propagation is precisely what turns an otherwise straightforward reliability assessment into a combinatorially complex problem, a complexity that we later make rigorous through an NP-hardness result. 
In this sense, causal failures add complexity to the analysis of network vulnerability.

In this work, we mainly consider the impact of cascading causal failures on the paths between a set of source and terminal nodes, with contributions as follows:
\begin{itemize}
    \item We formally define cascading causal failures. If the nodes in $V_1\subseteq V$ ($V$ is the set of network nodes) fail, all the nodes in $V_2\subseteq V$ will fail. 
    The nodes in $V_1$ do not need to be adjacent to the nodes in $V_2$.
    Furthermore, if $V_3 \subseteq V_1\cup V_2$ and the failure of the nodes in $V_3$ causes the failure of the nodes in $V_4$, then the failure of $V_1$ leads to cascading causal failures.
    \item We define the closure of a subset $U$ of the node set as the set of nodes whose failures are triggered by the failure of $U$.
    More specifically, the closure of $U\subseteq V$ based on the causality set is a set of failed nodes defined in the causalities and initially triggered by the failure of the nodes in $U$.
    \item We consider the impact of causal node failures on the shortest paths of a pair of nodes $s, t$ in the network $G$ and formulate a problem to find the maximum number of causal failures applied, restricting the paths $st$ to at most $\alpha \times d_G(s, t)$.
    \item We provide an NP-hardness proof of the formulated problem.
    \item We present an integer linear programming model that solves the proposed problem exactly.
    To scale to larger networks, we further design a column-generation-based approximation, together with a greedy algorithm as a fast alternative, and compare the three methods in terms of solution quality and running time.
\end{itemize}

This paper is organized as follows. Section \ref{S:2} introduces related work. Section \ref{S:3} shows the system model of this paper. We provide an NP-hardness proof of the proposed problem, the corresponding integer linear programming model, a column-generation-based approximation, and a greedy heuristic in Section \ref{S:4}. We evaluate the optimization model, the column generation algorithm and the greedy heuristic on several networks in Section \ref{S:6}. Finally, Section \ref{S:7} concludes this paper.

\section{Related Work}
\label{S:2}

The functionality of nodes is often correlated in networked systems, where the failure of some nodes can affect the functionality of other nodes, potentially causing failure in those nodes~\cite{Rinaldi2001}. 
In cyber-physical systems, one node functions in a certain state only when it belongs to the giant component of its network and has support from the nodes of the other networks \cite{yagan}. This scheme indicates that once there are node failures in one network, the other network may also have failed nodes due to interdependence (e.g., cascading failures)~\cite{Buldyrev2010}. Zhang et al. \cite{Zuyuan} analyzed reliability-based cascading failures in cyber-physical systems and provided an optimal inter-edge allocation strategy based on reliability.
In a power transmission system, Lin et al. \cite{LIN2012954} measured the impact of correlated failures on the reliability of a power network, where correlated failures refer to simultaneous failures of physical lines at each edge that connects two nodes with a given probability. 
In a computer network, Lin et al. \cite{LIN201332} presented correlated failures of transmission lines and routers and proposed an algorithm to assess the impact of correlated failures on the reliability of network flow.
More recently, Zhang et al. \cite{zzy_causal} proposed the causal failure model to address the dependence of nodes in a networked system and considered the robustness and vulnerability of a network based on the size of the giant component and the number of connected components, respectively.

The identification of the shortest paths has applications in a variety of engineering systems \cite{Cancela, Olivera, sed}. 
For a power system, Biswas et al. \cite{Biswas} proposed a new approach based on shortest paths to analyze whether a contingency creates a saturated cutoff in a power system when considering the problem of identifying transmission lines with limited capacity if multiple disturbances occur in the system.
For the singularity-induced voltage instability in a power network, Wu et al. \cite{WU2021107127} searched for the shortest path to its boundary. The edge length of the path is used to formulate the problem and is converted into a control framework to solve the shortest path manifold.
Ghasemkhani et al. \cite{Ghasemkhani} noted the importance of designing appropriate management, measurement, and communication infrastructures for a wide-area measurement system, utilizing dynamic multiobjective shortest-path programming for optimal communication infrastructure design.
For the virtualization of network functions, to find a suitable service path that optimizes a specific target, Sasabe and Hara \cite{Sasabe} formulated a capacitated integer linear programming model based on the shortest path tour problem for service chaining and function placement, and showed that the performance of the proposed model is better than the current state model.

Constraints on the length of the path require that the distance between any pair of nodes be at most a given integer $D$.
The diameter of a graph is the greatest distance between any pair of nodes, therefore a constraint on path length is referred to as a diameter constraint.
A diameter constraint has been widely used in many network analyses, especially when evaluating network reliability~\cite{Cancela2004}. 
Cancela et al. \cite{Petingi} introduced a diameter constraint on source-to-terminal reliability in a communication network, where each link fails with a given probability. 
They integrated the constraint with a factoring algorithm \cite{sat} and illustrated the gain in computation in terms of topological reduction when computing reliability.
In the same spirit of topological reduction, Xiong et al. \cite{Xiong} proposed a polynomial-time algorithm that detects edges irrelevant to diameter-constrained reliability, so that they can be removed without altering network reliability.
In a transportation network, which is often modeled as a flow network, Zhang and Shao \cite{zzy} designed an approximation algorithm for multistate two-terminal reliability using a diameter constraint of $s$--$t$ paths as a tuning parameter.
Moreover, due to signal attenuation problems \cite{Mutafungwa:01}, wavelength division multiplexing networks often require a diameter constraint. Similarly, a hop constraint is also useful in virtual topology and routing design \cite{Skorin-Kapov08, Leesutthipornchai, LEESUTTHIPORNCHAI20102246}.

\begin{table}[!t]
\centering
\caption{Contributions of related work to the modeling aspects considered in
this paper.}
\label{tab:related}
\footnotesize
\setlength{\tabcolsep}{3pt}
\begin{tabular}{@{}lccccc@{}}
\hline
References & \makecell{Network \\ failures} & Paths &
\makecell{Length \\ constraints} & \makecell{Causal \\ failures} &
\makecell{Column \\ generation} \\
\hline
\cite{yagan,Zuyuan} & \checkmark & & & & \\
\cite{LIN2012954,LIN201332} & \checkmark & \checkmark & & & \\
\cite{zzy_causal} & \checkmark & & & \checkmark & \\
\cite{Cancela,Olivera,sed} & & \checkmark & & & \\
\cite{Petingi,Cancela2004,Xiong,zzy} & \checkmark & \checkmark & \checkmark & & \\
\cite{Skorin-Kapov08,Leesutthipornchai,LEESUTTHIPORNCHAI20102246} & & \checkmark & \checkmark & & \\
\cite{Barnhart1998,Lubbecke2005,Ahmadi2025RCSP} & & \checkmark & & & \checkmark \\
\textbf{This paper} & \checkmark & \checkmark & \checkmark & \checkmark & \checkmark \\
\hline
\end{tabular}
\end{table}

A third thread of related work, methodological rather than topical, deals with column generation for path-based integer programs.
Column generation is a standard tool for integer programs with an exponential number of variables.
A restricted master problem is iteratively augmented with new columns priced by an auxiliary subproblem.
When embedded in a branch-and-bound scheme, it becomes
branch-and-price~\cite{Barnhart1998,Lubbecke2005}.
For path-based formulations, the pricing subproblem is typically a resource-constrained shortest-path problem.
Edge costs come from the dual variables of the master, and resource consumption is bounded along the path.
Label-setting algorithms and bidirectional variants are the standard
approach.
We use the recent algorithm for label-setting $A^{*}$ of Ahmadi et al.~\cite{Ahmadi2025RCSP} as the exact level of our pricing cascade.
To our knowledge, column generation has not previously been applied to length-constrained shortest-path problems under cascading causal failures.

\section{Formalism of the System Model}
\label{S:3}

\begin{figure}[!t]
    \centering
    \includegraphics[width=0.6\textwidth]{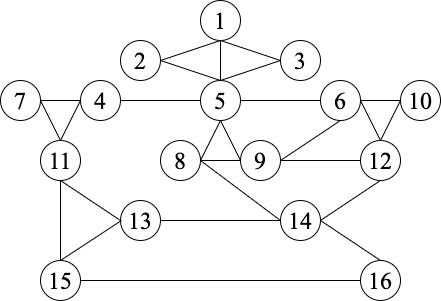}
    \caption{An example of a power network with 16 nodes and 26 edges.}
    \label{An example of a network}
\end{figure}

In this work, we mathematically model a network as an undirected graph $G = (V, E)$, where $V$ is the set of nodes with $|V|=n$, and $E$ is the set of edges with $|E| = m$.
Fig. \ref{An example of a network} depicts an example network with 16 nodes and 26 edges \cite{zzy_causal}.

We consider the causal dependence on the nodes \cite{zzy_causal}:
\begin{definition}
\label{causal_failures}
Given $V_1, V_2 \subset V$, if failures of all nodes in $V_1$ trigger failures of all nodes in $V_2$, we say that this is a causal failure from $V_1$ to $V_2$, denoted as $C: V_1\Rightarrow V_2$.
\end{definition}

Def. \ref{causal_failures} describes the failure of some nodes as a result of the failure of some other nodes and does not require that the nodes in $V_1$ and $V_2$ be adjacent. This indicates that Def. \ref{causal_failures} provides a general case for all kinds of node failure in a network.

\begin{figure}[!t]
    \centering
    \includegraphics[width=0.4\textwidth]{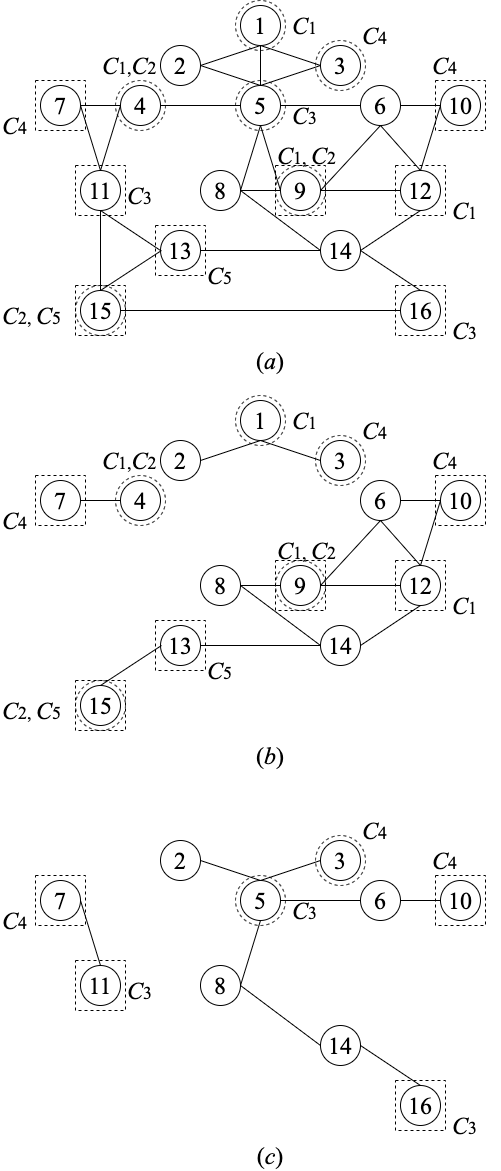}
    \caption{Remaining networks after applying different causalities.}
    \label{Causal_failures_Network}
\end{figure}

\begin{example}
\label{ex:1}
In Fig. \ref{Causal_failures_Network}($a$), we use dotted circles and dotted squares to denote the nodes in $V_1$ and $V_2$ of a causality $C$, respectively. 
The given causalities are $C_{1}: \{1, 4\}\Rightarrow \{9, 12\}$, $C_2: \{4, 9\}\Rightarrow \{15\}$, $C_{3}: \{5\}\Rightarrow \{11, 16\}$, $C_{4}: \{3\}\Rightarrow \{7, 10\}$, $C_5: \{15\}\Rightarrow\{13\}$.

If we apply $C_{3}$ to the network, the failure of node 5 causes the failure of nodes 11 and 16. 
Fig. \ref{Causal_failures_Network}($b$) presents the remaining network with three connected components (i.e., $\{1, 2, 3\}$, $\{4, 7\}$, and $\{6, 8, 9, 10, 12, 13, 14, 15\}$). 

If $C_1$ is applied to the network, then nodes $\{1, 4, 9, 12\}$ are removed from the network.
As $\{4, 9\}$ is on the left side of $C_2$, node 15 fails.
Furthermore, failure of node 15 causes failure of node 13, according to $C_5$.
Therefore, we observe a cascade of causal failures $C_1, C_2, C_5$, and Fig. \ref{Causal_failures_Network}($c$) shows the remaining network consisting of two connected components (that is, $\{2, 3, 5, 6, 8, 10, 14, 16\}$ and $\{7, 11\}$) after applying $C_1, C_2, C_5$.
\end{example}

From Example \ref{ex:1}, we observe cascading failures in the causal model.
Therefore, we formally define cascading causal failures.

\begin{definition}
\label{Causal_Cascading_Failures}
Given a network $G = (V, E)$ and a causality set $\mathcal{C}$, if the application of a causality $C\in \mathcal{C}$ triggers the subsequent application of other causalities, we say that cascading causal failures are applied to the network.
\end{definition}

Def. \ref{Causal_Cascading_Failures} generalizes causal failures to the cascading case if the information on the dependency between nodes is known.
Naturally, given a group of failed nodes, it is necessary to obtain all failed nodes triggered by the failures of the group of nodes based on causal failures.
Therefore, we define the closure of a group of nodes.

\begin{definition}
\label{Closure}
Given a set of causalities $\mathcal{C}$, the closure of a set of nodes $U$ with respect to $\mathcal{C}$, denoted by $U^{+}_{\mathcal{C}}$, 
is the set of all nodes (including those in $U$) whose failures are triggered by causalities in $\mathcal{C}$ after removal of nodes in $U$. 
The cascading set $\mathcal{C}(U)$ is the set of causalities in $\mathcal{C}$ that are automatically triggered (in a cascading fashion) as a result of the failure of the nodes in $U$.  
\end{definition} 

The closure plays an important role in cascading causal failures assuming the time delays for cascading effects to be negligible while studying the performance of the network. 
An equivalent algorithmic definition of $U^+_{\mathcal{C}}$ and $\mathcal{C}(U)$ is the result of Algorithm 1, the computational complexity of which is $O(h^2)$. 

\begin{algorithm}[!t]
\caption{Closure}
\label{closure}
\textbf{Input}: graph $G = (V, E)$, $U \subseteq V$, a causality set $\mathcal{C}$;\\
\textbf{Output:} $U^+_{\mathcal{C}}$, $\mathcal{C}(U)$. 
\begin{algorithmic}[1]
\State  $U' \leftarrow U$, $S \leftarrow \emptyset$, $\mathcal{D} \leftarrow \mathcal{C}$ ; 
\While{$\exists C_i \in \mathcal{D}$ such that $V_{1,i} \subseteq U'$}
\State $U' \leftarrow U' \cup V_{2,i},  S \leftarrow S \cup \{C_i\},  
\mathcal{D} \leftarrow \mathcal{D} \setminus \{C_i\}$
\EndWhile
\State \textbf{return} $U'$ and $S$. 
\end{algorithmic} 
\end{algorithm}

In Example \ref{ex:1}, suppose $U=\{1, 4\}$, then according to the causality set, causalities $C_1, C_2, C_5$ would be applied automatically if the nodes in $U$ fail.
Hence, given that $U=\{1,4\}$, we have $\mathcal{C}(U)=\{C_1, C_2, C_5\}$ and $U^{+}_{\mathcal{C}}=\{1, 4, 9, 12, 13, 15\}$ running Algorithm 1.

Given a source node $s$ and a terminal node $t$, it is well known that a node failure may increase the distance between them \cite{zzy,Albert2000}. 
More specifically, the more node failures the network has, the more likely that the distance between $s$ and $t$ increases.

\begin{example}
\label{ex:2}
In Fig. \ref{Causal_failures_Network}, let $s = 2$ and $t = 14$, then $d_{G}(2, 14) = 3$. If causality $C_1$ is applied to $G$, then nodes 1, 4, 9, 12, 15, and 13 will be removed from $G$, but the distance between 2 and 14 is still 3. However, if we apply the causal failure $C_3$, then node $2$ is not reachable from node $14$ because node $5$ is removed. 
\end{example}

Example \ref{ex:2} presents the difference in causalities with respect to their impact on connectivity. 
As alluded to in Section 1, the lengthening of the shortest paths can adversely affect certain networks (e.g., the distance between nodes affects the efficiency of power transmission in the power network~\cite{Latora2001,Crucitti2004}). Therefore, it is interesting to investigate how the number of causal failures applied affects the shortest path from $s$ to $t$. 
Leaving $D_0 = d_G(s, t)$, we have the following problem.

\begin{problem}
\label{Pro:1}
Given a network $G=(V, E)$, causality set $\mathcal{C} = \{C_{i}: V_{1,i}\Rightarrow V_{2,i}, i = 1, 2, \cdots, h\}$ (containing possible cascading failures), $s, t\in V$, and a constant $\alpha$, find a {\bf maximum} subset of causalities whose application maintains the distance between $s$ and $t$ to be at most $\alpha \times D_0$.
\end{problem}

Problem \ref{Pro:1} investigates the impact of cascading causal failures on the distance between a given pair of nodes. 
This is a good metric to measure the connectivity of a pair of nodes in terms of the maximum number of causalities they can withstand.
Let $D = \alpha \times D_0$, and we shall use $D$ and $\alpha$ interchangeably below.

\section{Problem Analysis and Algorithms}
\label{S:4}

Problem \ref{Pro:1} is also applicable to edge-weighted networks because it focuses on the length of the $s$--$t$ paths.
To prove its hardness, we start with the case where the edges are not weighted and there are no cascading failures and then extend the unweighted case to the weighted case.

\subsection{NP-Hardness Proof}

We consider the decision version of Problem \ref{Pro:1} with unweighted edges and no cascading failures, denoted as \textit{VfCD}: given a network $G = (V, E)$, causalities $\mathcal{C} = \{C_{i}: V_{1,i}\Rightarrow V_{2,i}, i = 1, 2, \ldots, h\}$, $s, t\in V$, $D_0 = d_G(s, t)$, a constant $\alpha \geq 1$, and an integer $f$, let $D = \floor*{\alpha\times D_0}$; does there exist $\mathcal{C}^{\prime} \subseteq \mathcal{C}$ with $|\mathcal{C}^{\prime}|\geq f$ such that applying the causalities in $\mathcal{C}^{\prime}$ --- that is, removing the node set $\bigcup_{C\in\mathcal{C}^{\prime}}(V_{1,C}\cup V_{2,C})$ from $G$ --- maintains $d_{G}(s, t)\leq D$?

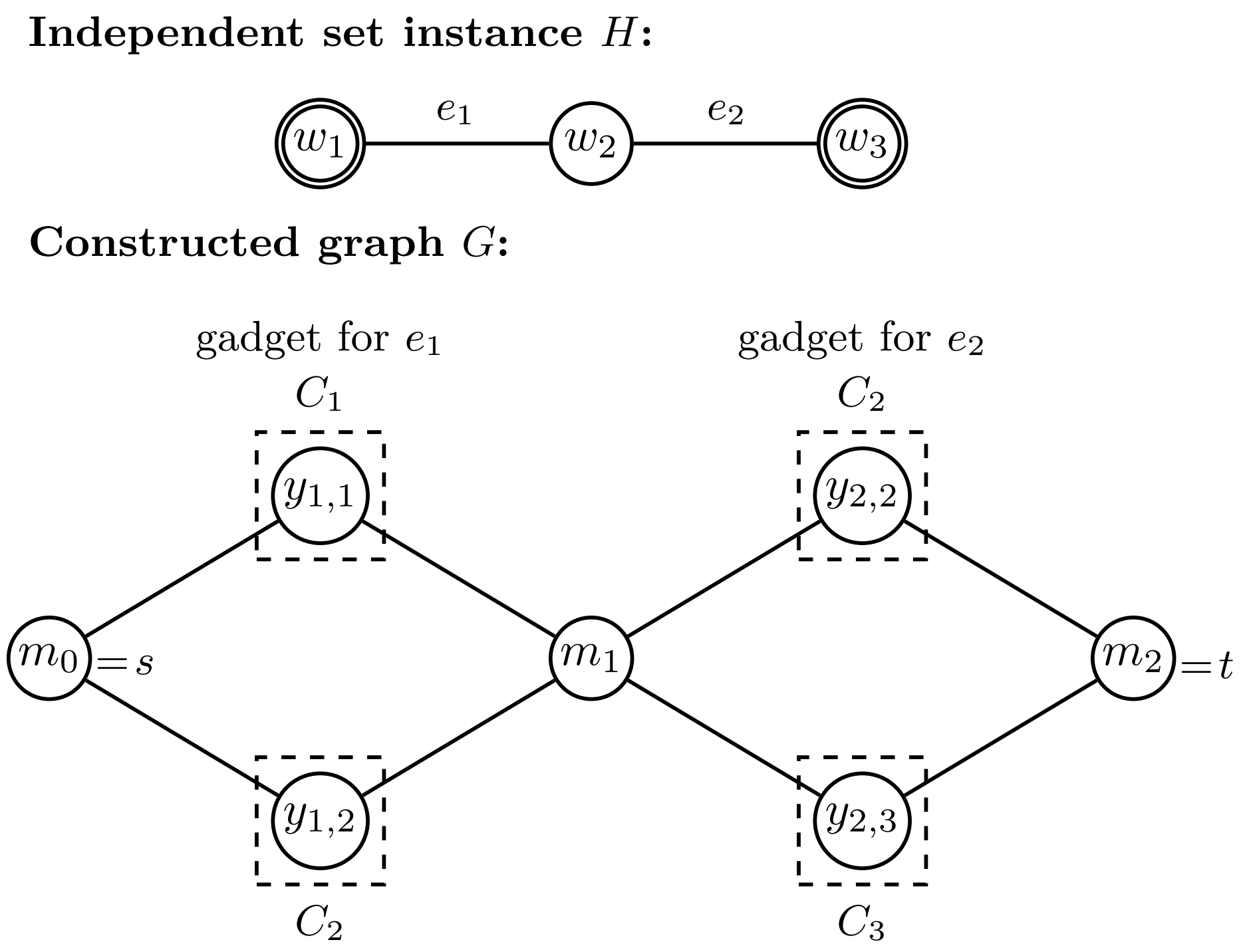
\begin{figure}[!t]
    \centering
    \resizebox{0.75\textwidth}{!}{%
    \begin{tikzpicture}[
        vtx/.style={draw,circle,inner sep=1pt,minimum size=17pt,font=\small},
        caus/.style={draw,dashed,inner sep=3pt},
        thick]
        \node[font=\footnotesize\bfseries,anchor=west] at (-0.3,4.6) {Independent set instance $H$:};
        \node[vtx,double] (w1) at (2,3.8) {$w_1$};
        \node[vtx]        (w2) at (4,3.8) {$w_2$};
        \node[vtx,double] (w3) at (6,3.8) {$w_3$};
        \draw (w1) -- node[above,font=\footnotesize] {$e_1$} (w2)
                  -- node[above,font=\footnotesize] {$e_2$} (w3);
        \node[font=\footnotesize\bfseries,anchor=west] at (-0.3,3.05) {Constructed graph $G$:};
        \node[font=\footnotesize] at (2,2.35) {gadget for $e_1$};
        \node[font=\footnotesize] at (6,2.35) {gadget for $e_2$};
        \node[vtx] (m0) at (0,0) {$m_0$};
        \node[vtx] (m1) at (4,0) {$m_1$};
        \node[vtx] (m2) at (8,0) {$m_2$};
        \node[font=\footnotesize,anchor=west] at (0.25,-0.05) {$\!=\!s$};
        \node[font=\footnotesize,anchor=west] at (8.25,-0.05) {$\!=\!t$};
        \node[vtx] (y11) at (2,1.2)  {$y_{1,1}$};
        \node[vtx] (y12) at (2,-1.2) {$y_{1,2}$};
        \node[vtx] (y22) at (6,1.2)  {$y_{2,2}$};
        \node[vtx] (y23) at (6,-1.2) {$y_{2,3}$};
        \draw (m0) -- (y11) -- (m1);
        \draw (m0) -- (y12) -- (m1);
        \draw (m1) -- (y22) -- (m2);
        \draw (m1) -- (y23) -- (m2);
        \node[caus,fit=(y11),label={[font=\footnotesize]above:$C_1$}]{};
        \node[caus,fit=(y22),label={[font=\footnotesize]above:$C_2$}]{};
        \node[caus,fit=(y12),label={[font=\footnotesize]below:$C_2$}]{};
        \node[caus,fit=(y23),label={[font=\footnotesize]below:$C_3$}]{};
    \end{tikzpicture}}
    \caption{Construction in the proof of Theorem~\ref{no_cascade_proof}, illustrated on the \textsc{Independent Set} instance $H$ with $W=\{w_1,w_2,w_3\}$ and edges $e_1=\{w_1,w_2\}$, $e_2=\{w_2,w_3\}$ ($M=2$); the independent set $\{w_1,w_3\}$ is drawn with double circles. Each edge $e_r$ of $H$ becomes a gadget of two parallel length-two routes, placed in series on the $s$--$t$ spine $m_0,m_1,m_2$, so $D_0=d_G(s,t)=2M=4$; with $\alpha=1$ this gives $D=4$ (any $\alpha\geq 1$ gives the same verdicts, as $D\geq 2M$). As in Fig.~\ref{Causal_failures_Network}, a dashed box marks each node removed by a causality and is labelled with that causality: $C_i$ removes the route node $y_{r,i}$ of every gadget incident to $w_i$ (the isolated trigger nodes $g_i$ are omitted, as they affect no distance). A subset of causalities keeps an $s$--$t$ path of length at most $D$ if and only if the corresponding nodes form an independent set of $H$. For example, $\{C_1,C_3\}$ removes $y_{1,1},y_{2,3}$ and leaves the path $m_0\,y_{1,2}\,m_1\,y_{2,2}\,m_2$ of length $4$, whereas $\{C_1,C_2\}$ removes both $y_{1,1}$ and $y_{1,2}$, severing $m_0$ from $m_1$ and disconnecting $s$ from $t$.}
    \label{fig:thm1}
\end{figure}

To prove that this version is NP-complete, we reduce from the \textsc{Independent Set} problem~\cite{Karp1972}, which is NP-complete: given an undirected graph $H=(W,F)$ and an integer $k$, decide whether $H$ contains an independent set of size at least $k$ (i.e.,\ a set of pairwise non-adjacent nodes). We may assume $|F|\geq 1$, since the problem is trivial otherwise, and that $H$ has no isolated node, since an isolated node belongs to every maximum independent set and can be removed after decreasing $k$ by one.
\begin{theorem}
\label{no_cascade_proof}
VfCD is NP-complete.
\end{theorem}
\begin{proof}
\emph{Membership.} Given a certificate $\mathcal{C}'\subseteq\mathcal{C}$, one checks $|\mathcal{C}'|\geq f$, removes the node set $\bigcup_{C\in\mathcal{C}'}(V_{1,C}\cup V_{2,C})$ from $G$, and computes $d_{G}(s,t)$ in the resulting graph by breadth-first search, accepting if and only if $d_{G}(s,t)\leq D$. Every step is polynomial, so \textit{VfCD}$\,\in\,$NP.

\emph{Hardness.} Let $(H=(W,F),k)$ be an instance of \textsc{Independent Set}, with $W=\{w_1,\dots,w_q\}$ and $F=\{e_1,\dots,e_M\}$, $M\geq 1$. We build a \textit{VfCD} instance as follows (Fig.~\ref{fig:thm1} gives an example).
\begin{enumerate}
    \item \emph{Spine.} Create nodes $m_0,m_1,\dots,m_M$ and set $s=m_0$ and $t=m_M$.
    \item \emph{Edge.} For each edge $e_r=\{w_i,w_j\}$ ($r=1,\dots,M$) create two nodes $y_{r,i},y_{r,j}$ and the four edges $(m_{r-1},y_{r,i})$, $(y_{r,i},m_r)$, $(m_{r-1},y_{r,j})$, $(y_{r,j},m_r)$. This places two internally disjoint routes of length two between $m_{r-1}$ and $m_r$, and the gadgets are chained in series along the spine.
    \item \emph{Causalities.} For each node $w_i\in W$ create a private node $g_i$, isolated in $G$, and one causality
    \[ C_i:\ \{g_i\}\Rightarrow\{\,y_{r,i}\ :\ e_r \text{ is incident to } w_i\,\}. \]
    Because $H$ has no isolated node, every effect set is nonempty. The trigger node $g_i$ is private to $C_i$ and lies in the effect set of no causality, so applying any subset of causalities never removes the trigger of another one; hence the instance has no cascading failures and the closure of $V(\mathcal{C}')$ equals $\bigcup_{C\in\mathcal{C}'}(V_{1,C}\cup V_{2,C})$.
    \item Set $f=k$. All edges have unit weight, hence $D_0=d_{G}(s,t)=2M$ and $D=\floor*{\alpha\times 2M}\geq 2M$ for every constant $\alpha\geq 1$.
\end{enumerate}
The construction has $3M+q+1=O(M+q)$ nodes, $4M$ edges and $q$ causalities, and is produced in polynomial time.

Since there is exactly one causality $C_i$ per node $w_i$, the map $\mathcal{C}'\mapsto S=\{w_i:C_i\in\mathcal{C}'\}$ is a bijection between subsets of $\mathcal{C}$ and subsets of $W$ with $|\mathcal{C}'|=|S|$. Applying $\mathcal{C}'$ removes the isolated nodes $\{g_i:w_i\in S\}$, which affect no distance, together with every $y_{r,i}$ such that $w_i\in S$. The spine nodes $m_0,\dots,m_M$ belong to no causality and are never removed. Each $y_{r,i}$ has degree two, its only neighbours being $m_{r-1}$ and $m_r$, and the graph has no further edges; hence for $0<r<M$ the node $m_r$ separates $\{m_0,\dots,m_{r-1}\}$ from $\{m_{r+1},\dots,m_M\}$, so every $s$--$t$ path traverses $m_0\to m_1\to\cdots\to m_M$ in order, using one route of each gadget (for $M=1$ this is immediate, as $s=m_0$ and $t=m_1$ are the two ends of the single gadget). Deleting nodes cannot create new connections, so this remains true after applying $\mathcal{C}'$. We show that
\begin{gather*}
\text{applying } \mathcal{C}' \text{ keeps } d_{G}(s,t)\leq D\\
\Longleftrightarrow\quad
S \text{ is an independent set of } H.
\end{gather*}

$(\Leftarrow)$ If $S$ is independent, then for every edge $e_r=\{w_i,w_j\}$ at most one of $w_i,w_j$ lies in $S$, so at least one of $y_{r,i},y_{r,j}$ survives and gadget $r$ still offers a length-two route from $m_{r-1}$ to $m_r$. Concatenating these routes for $r=1,\dots,M$ yields an $s$--$t$ path of length $2M\leq D$; hence $d_{G}(s,t)\leq D$.

$(\Rightarrow)$ If $S$ is not independent, some edge $e_r=\{w_i,w_j\}$ has both endpoints in $S$, so both $y_{r,i}$ and $y_{r,j}$ are removed. As these are the only internal nodes of the two routes joining the consecutive cut nodes $m_{r-1}$ and $m_r$, no path links $m_{r-1}$ to $m_r$; thus $s$ and $t$ are disconnected and $d_{G}(s,t)=\infty>D$.

Therefore the largest applicable subset of causalities has size exactly the maximum independent set of $H$, and a feasible $\mathcal{C}'$ with $|\mathcal{C}'|\geq f$ exists if and only if $H$ has an independent set of size at least $k=f$. Since \textsc{Independent Set} is NP-complete, \textit{VfCD} is NP-hard; together with membership, \textit{VfCD} is NP-complete.
\end{proof}

Theorem \ref{no_cascade_proof} establishes the hardness of Problem \ref{Pro:1} with unweighted edges. The unweighted case is exactly the special case of unit edge weights, hence a restriction of the edge-weighted problem, and we immediately obtain:
\begin{corollary}
VfCD with weighted edges is NP-complete.
\end{corollary}

Finally, Theorem~\ref{no_cascade_proof} and Corollary~1 concern the version of Problem~\ref{Pro:1} without cascading failures, which is the special case in which every causality is triggered in isolation; the general problem therefore inherits its hardness. Membership in NP still holds, since the closure of a node set is computed in polynomial time by Algorithm~1. Therefore, we have the following:
\begin{corollary}
VfCD with weighted edges and cascading failures is NP-complete.
\end{corollary}

\subsection{Optimization Model}

We propose a mixed integer linear programming (MILP) formulation to find the maximum number of causal failures that result in the shortest path of at most $D=\alpha\times D_0$. We consider the directed version $DG = (V, A)$ of the undirected graph $G = (V, E)$, where $A$ is generated by replacing each edge in $E$ with two directed edges in opposite directions. Given the causal failures set $\mathcal{C}$ and the length constraint $D$, we have the formulation in Eqs.~\eqref{eq:milp-obj}--\eqref{eq:milp-domy}.

We first list the decision variables and parameters used in the formulation.

\begin{center}
\footnotesize
    \begin{tabular}{|c|p{0.6\textwidth}|}
    \hline
     Decision variables & Definition of decision variables \\ \hline
     $z_C$   & A binary variable indicating whether causality $C$ is applied \\ \hline
     $x_{ij}$ & A binary variable indicating the flow of edge $(i,j)$\\ \hline
     $y_{i}$ & A binary variable indicating the integrity of node $i$\\ \hline
    \end{tabular}
\end{center}

\begin{center}
\footnotesize
    \begin{tabular}{|c|p{0.6\textwidth}|}
    \hline
     Parameters & Definition of parameters \\ \hline
     $w_{ij}$   & The weight of edge $(i,j)$ \\ \hline
     $i, j$ & Nodes $i, j$\\ \hline
     $D$ & The length constraint \\ \hline
     $C$ & A causality in $\mathcal{C}$\\ \hline
    \end{tabular}
\end{center}

A central modeling choice is to represent the integrity of each node using a dedicated binary variable, $y_i$, kept separate from the routing variables $x_{ij}$. The routing variables alone cannot express the removal of a node: an edge that does not carry flow ($x_{ij} = 0$) may be unused by the chosen path or incident to a node that a causality has removed. 
Conflating these two cases would allow a healthy node that simply lies off the path to be read as failed, spuriously triggering its causalities. 
The variable $y_i$ removes this ambiguity by explicitly encoding the state of the node $i$: $y_i = 1$ when the node is intact and $y_i = 0$ when it has been removed by an applied causality. 
The routing is then admitted through a node only when the node is intact.
A causality is triggered only when its trigger nodes are genuinely removed.

The MILP formulation is as follows.
\begin{subequations}\label{prob:milp}
\begin{align}
\max \quad & \sum_{C\in \mathcal{C}}z_C \label{eq:milp-obj} \\
\text{s.t.} \quad
& \sum_{ij\in A}w_{ij}x_{ij} \leq D \label{eq:milp-dist} \\
& \sum_{j\in V}x_{ij} - \sum_{j\in V}x_{ji} =
  \begin{cases}
    1, & i = s\\
    -1, & i = t\\
    0, & \text{otherwise}
  \end{cases} \quad \forall i \in V \label{eq:milp-flow} \\
& \sum_{j \in V}x_{ij} \leq 1 \quad \forall i\in V\setminus \{s, t\} \label{eq:milp-deg} \\
& y_i \leq 1 - z_C \quad \forall C\in \mathcal{C}, \; i \in V_{1,C} \cup V_{2,C} \label{eq:milp-eff} \\
& \sum_{j\in V}x_{ij} \leq y_i \quad \forall i\in V \label{eq:milp-fout} \\
& \sum_{j\in V}x_{ji} \leq y_i \quad \forall i\in V \label{eq:milp-fin} \\
& z_C \geq 1 - \sum_{i\in V_{1, C}}y_i \quad \forall C\in \mathcal{C} \label{eq:milp-trig} \\
& x_{ij} \in \{0, 1\} \quad \forall (i, j)\in A \label{eq:milp-domx} \\
& z_C \in \{0, 1\} \quad \forall C\in \mathcal{C} \label{eq:milp-domz} \\
& y_i \in \{0, 1\} \quad \forall i\in V \label{eq:milp-domy}
\end{align}
\end{subequations}

This formulation aims to transmit one unit of flow from source $s$ to terminal $t$. 
Eq.~\eqref{eq:milp-obj} is the objective function, which is to maximize the number of causalities that could be applied.
We require the length of the path to be at most $D$, according to Eq.~\eqref{eq:milp-dist}.
In Eq.~\eqref{eq:milp-flow}, for the source $s$, the sum of the outgoing flow is 1 higher than the incoming flow, while the total incoming flow for $t$ is 1 higher than its total outgoing flow. Other nodes follow the flow conservation constraint. 
Because only one unit of flow is transmitted in the network, the total outgoing flow of each node $i\in V\setminus\{s, t\}$ is at most 1, as ensured by Eq.~\eqref{eq:milp-deg}.
Consequently, the state of each edge would be $0$ or $1$, and the edges with flow $1$ definitely form an $s$--$t$ path. 
Regarding causality $C$, the binary variable $z_C = 1$ implies that causality $C$ is applied to the network.
Eqs.~\eqref{eq:milp-eff}--\eqref{eq:milp-trig} tie the routing to the integrity of the nodes and govern the application of cascading causal failures.
By Eq.~\eqref{eq:milp-eff}, if a causality $C$ is applied ($z_C = 1$), then every node involved in $C$ (i.e., in $V_{1,C}\cup V_{2,C}$) is removed ($y_i = 0$).
Eqs.~(1f) and (1g) state that a node can carry outgoing or incoming flow only if it is intact.
No flow may cross a node with $y_i = 0$.
An intact node ($y_i = 1$) may still carry flow or none if the path does not use it.
Finally, Eq.~\eqref{eq:milp-trig} enforces the cascading effect: if all trigger nodes of causality $C$ are removed, that is, $\sum_{i\in V_{1, C}}y_i = 0$, then $z_C$ is forced to $1$ and $C$ is applied automatically.
The decoupling of $y_i$ from the routing variables $x_{ij}$ is precisely what makes Eq.~\eqref{eq:milp-trig} correct: a node that merely lies off the path keeps $y_i = 1$, so a causality is triggered only by an actual removal of its trigger nodes and never by the routing choice.

\subsection{Column Generation-Based Approximation Algorithm}

We now present an alternative, path-based formulation that is amenable to column generation. Let $\mathcal{P}$ be the set of $s$--$t$ paths of length at most $D$, and for $P\in\mathcal{P}$ let $\ell_{P} = \sum_{(i,j)\in P} w_{ij}$ be its length. Because $|\mathcal{P}|$ can grow exponentially with the size of the network, we do not enumerate it explicitly: we keep a small subset $\mathcal{P}' \subseteq \mathcal{P}$ and expand it on demand following the principle of column generation. We introduce a variable $x_{P}$ for the weight of the path $P$ and reuse the node-integrity variables $y_i$ and the causality variables $z_C$ of the previous model.

The linear relaxation of the restricted master problem (RMP) over $\mathcal{P}'$ is as follows, where the dual variable of each constraint is reported on the right.
\begin{subequations}\label{prob:cg-master}
\begin{align}
\max \quad & \sum_{C\in \mathcal{C}}z_C \label{eq:cg-obj} \\
\text{s.t.} \quad
& \sum_{P\in \mathcal{P}'} x_{P} = 1 \quad [\lambda \text{ free}] \label{eq:cg-conv} \\
& \sum_{P\in \mathcal{P}'} \ell_{P}\, x_{P} \leq D \quad [\mu \geq 0] \label{eq:cg-len} \\
& \sum_{P\in \mathcal{P}':\, i\in P} x_{P} \leq y_i \quad \forall i\in V \;\; [\gamma_i \geq 0] \label{eq:cg-cpl} \\
& y_i \leq 1 - z_C \quad \forall C\in \mathcal{C}, \; i \in V_{1,C} \cup V_{2,C} \label{eq:cg-eff} \\
& z_C \geq 1 - \sum_{i\in V_{1, C}} y_i \quad \forall C\in \mathcal{C} \label{eq:cg-trig} \\
& x_{P} \geq 0 \quad \forall P\in \mathcal{P}' \label{eq:cg-domx} \\
& y_i,\, z_C \in [0, 1] \label{eq:cg-domyz}
\end{align}
\end{subequations}

The first constraint selects one unit of $s$--$t$ path weight, and the second keeps its length within the budget $D$. The coupling constraint $\sum_{P:\, i\in P} x_{P} \leq y_i$ prevents any selected path from traversing a removed node ($y_i = 0$), while the last two constraints reproduce the effect and the trigger of cascading causal failures, exactly as in Eqs.~\eqref{eq:milp-eff} and~\eqref{eq:milp-trig}. Replacing $\mathcal{P}'$ with the complete set $\mathcal{P}$ yields the complete master problem.

The path variables $x_{P}$ occur only in the first three constraints, whose dual variables are $\lambda$ (free), $\mu \geq 0$ and $\gamma_i \geq 0$, respectively. The reduced cost of a path $P\in\mathcal{P}$ is therefore the case.
\begin{equation}\label{eq:cg-rc}
\bar{c}_{P} = \lambda - \mu\, \ell_{P} - \sum_{i\in P} \gamma_i .
\end{equation}
The current RMP solution is optimal for the entire master when $\bar{c}_{P} \leq 0$ for all $P\in\mathcal{P}$; otherwise, a path with $\bar{c}_{P} > 0$ is added as a new column. Since $\lambda$ is constant, generating such a path amounts to the pricing problem.
\begin{subequations}\label{prob:cg-pricing}
\begin{align}
\min_{P\in\mathcal{P}} \quad & \mu\, \ell_{P} + \sum_{i\in P} \gamma_i \label{eq:cg-pobj} \\
\text{s.t.} \quad & \ell_{P} \leq D \label{eq:cg-pcon}
\end{align}
\end{subequations}
That is, the above is a resource-constrained shortest-path problem with nonnegative node costs $\gamma_i$, edge costs proportional to $w_{ij}$, and a single length resource bounded by $D$. An improving column exists if and only if its optimal value is smaller than $\lambda$. 
We solve the pricing problem with a three-level cascade.
The first level is a beam-search heuristic that explores a bounded number of partial paths and returns immediately if it finds one with positive reduced cost.
The second level is a label-setting $A^{*}$ search for the
resource-constrained shortest-path problem~\cite{Ahmadi2025RCSP}, capped at a label budget that is set as a function of the network size.
If the cap is reached without certifying optimality, the search saves its label state and returns the best column found so far.
The third level is the same $A^{*}$ search without the cap, hot started from the saved state of the second level rather than restarted from scratch.
Each level passes control to the next only when no improving column has been found.
The cascade therefore returns a column as soon as a sufficiently good one exists and falls back to the exact pricer only when the heuristic levels
fail.

The generated columns are finally used to solve the path-based master in integer variables $x_{P}, y_i, z_C \in \{0, 1\}$, restricted to the pool $\mathcal{P}'$. This returns a feasible solution to the original problem, but not necessarily an optimal one: the columns in $\mathcal{P}'$ are priced out only with respect to the linear relaxation, so the integer optimum may require paths outside $\mathcal{P}'$. Closing this gap would require embedding the pricing within a branch-and-price scheme~\cite{Barnhart1998,Lubbecke2005}; instead, keeping the restricted integer master is what makes the overall method an approximation.

\subsection{Greedy Approximation Algorithm}

Given the NP-hardness of the formulated problem, we design a greedy strategy to solve the problem approximately in polynomial time. 
Let $\theta(C)$ be the distance between $s$ and $t$ if we apply the causality $C$ with its closure to $G$, then we sort the causalities based on the corresponding values of $\theta$ in ascending order.
In the sorted causalities, we add one causality and check the closure of the removed nodes and the distance between $s$ and $t$.
If the distance between $s$ and $t$ is less than or equal to $D$, we add the applied causalities to the result; otherwise, we ignore the present causality and continue to check the next one.
The details of the algorithm are listed below.

\begin{algorithm}[!t]
\caption{Greedy-Max-Length}
\label{alg:2}
\textbf{Input}: graph $G = (V, E)$, causality set $\mathcal{C}$, length constraint $D$, source node $s$ and terminal node $t$;\\
\textbf{Output:} a set of causalities $\mathcal{C}'$, closed under cascading, such that after applying all of them to $G$ there still exists an $s$--$t$ path of length at most $D$.
\begin{algorithmic}[1]
\State Sort the set $\mathcal{C}$ based on $\theta(C)$ for all $C\in \mathcal{C}$ and let $\theta(C_1)\leq \theta(C_2)\leq \cdots \leq \theta(C_h)$ and update $\mathcal{C}$ as $\{C_1, C_2, \cdots, C_h\}$;
\State $i \leftarrow 1$, $\mathcal{C}'\leftarrow \varnothing$;
\While{$i \leq h$}
\If{$C_i \notin \mathcal{C}'$}
\State Run Algorithm 1 to get the closure of $V(\mathcal{C}'\cup\{C_i\})$ (i.e., the causality set $S$ and the failed-node set $U$);
\If{$s, t$ are connected and $d_G(s, t)\leq D$ in $G-U$}
\State $\mathcal{C}' \leftarrow S$
\EndIf
\EndIf
\State $i \leftarrow i + 1$
\EndWhile \\
\Return $\mathcal{C}'$
\end{algorithmic} 
\end{algorithm}

\textbf{Complexity Analysis:} Let $n=|V|$, $m=|E|$ and $h=|\mathcal{C}|$. In Step~1, for each of the $h$ causalities Algorithm 2 computes its closure by Algorithm~1 in $O(h^2)$ time and one shortest-path computation for $d_G(s,t)$, and then sorts them. A single shortest-path computation costs $T_{\mathrm{sp}}=O(n+m)$ on unweighted graphs (breadth-first search) and $T_{\mathrm{sp}}=O(m+n\log n)$ on weighted graphs (Dijkstra's algorithm). Step~1 therefore costs $O\!\big(h(h^2+T_{\mathrm{sp}})\big)$, and the while loop runs at most $h$ iterations, each computing one closure in $O(h^2)$ and at most one shortest-path computation. Hence the overall time complexity of Algorithm 2 is $O\!\big(h^3 + h\,T_{\mathrm{sp}}\big)$, that is, $O\!\big(h^3 + h(n+m)\big)$ in the unweighted case and $O\!\big(h^3 + h(m+n\log n)\big)$ in the weighted case; in both cases it is polynomial.

\begin{theorem}
\label{thm:greedy-ratio}
Algorithm 2 returns a feasible solution, and for $h=|\mathcal{C}|\geq 2$, its worst-case approximation ratio is $h-1$.
\end{theorem}
\begin{proof}
Feasibility is immediate: a causality is accepted only when, after removing the failed-node set $U$ of the current closure, $s$ and $t$ stay at distance at most $D$.

Let $\mathrm{ALG}$ and $\mathrm{OPT}$ be the number of causalities returned by Algorithm 2 and by an optimal solution, respectively, and assume $\mathrm{OPT}\geq 1$ (otherwise both are $0$). Removing nodes can only increase the $s$--$t$ distance, so feasibility is monotone: if a set of causalities keeps $d_G(s,t)\leq D$, so does every subset. Two consequences follow. First, since $\mathrm{OPT}\geq 1$, at least one causality is applicable on its own, so the smallest-$\theta$ causality $C_1$ satisfies $\theta(C_1)\leq D$ and is accepted; hence $\mathrm{ALG}\geq 1$. Second, if all $h$ causalities are jointly applicable ($\mathrm{OPT}=h$), then by monotonicity each one is still applicable when it is examined, so Algorithm 2 accepts all of them and $\mathrm{ALG}=h$; equivalently, $\mathrm{ALG}=1$ forces $\mathrm{OPT}\leq h-1$. Combining, $\mathrm{OPT}/\mathrm{ALG}\leq h-1$ when $\mathrm{ALG}=1$, and $\mathrm{OPT}/\mathrm{ALG}\leq h/2\leq h-1$ when $\mathrm{ALG}\geq 2$ (using $\mathrm{OPT}\leq h$ and $h\geq 2$); thus $\mathrm{OPT}/\mathrm{ALG}\leq h-1$ in every case.

The ratio is tight. Let the only two $s$--$t$ paths be the internally disjoint $P_0=s\,a\,t$ and $P_1=s\,b\,c\,t$, and set $D=3$ (so $D_0=2$ and $\alpha=1.5$). Define $h$ causalities, each with its own private trigger node: $C_1$ removes $b$, while $C_2,\dots,C_h$ each remove $a$. Then $\theta(C_1)=2$ (only $P_0$ survives) and $\theta(C_j)=3$ for $j\geq 2$ (only $P_1$ survives), so Algorithm 2 examines $C_1$ first and accepts it; afterwards every $C_j$ would additionally remove $a$, leaving both $P_0$ (no $a$) and $P_1$ (no $b$) broken, so all are rejected and $\mathrm{ALG}=1$. The optimum instead applies $C_2,\dots,C_h$, which remove only $a$ and keep $P_1$, so $\mathrm{OPT}=h-1$.
\end{proof}

The factor $h-1$ is weak, but no efficient algorithm can do essentially better: the problem is as hard to approximate as \textsc{Maximum Independent Set}.

\begin{theorem}
\label{thm:inapprox}
Unless $\mathrm{P}=\mathrm{NP}$, for every constant $\varepsilon>0$, Problem~\ref{Pro:1} admits no polynomial-time $h^{1-\varepsilon}$-approximation algorithm, where $h=|\mathcal{C}|$.
\end{theorem}
\begin{proof}
The reduction of Theorem~\ref{no_cascade_proof} maps an \textsc{Independent Set} instance $H=(W,F)$ (without isolated nodes) to an instance of Problem~\ref{Pro:1} with exactly one causality per node, so $h=|W|$. As established there, applying a subset $\mathcal{C}'$ keeps $d_G(s,t)\leq D$ if and only if the corresponding node set $S$ is independent in $H$, and the map $\mathcal{C}'\mapsto S$ satisfies $|\mathcal{C}'|=|S|$. Hence the maximum number of applicable causalities equals the independence number $\alpha(H)$, and any feasible solution of value $v$ is, in polynomial time, an independent set of size $v$. A polynomial-time $\rho$-approximation for Problem~\ref{Pro:1} would therefore yield an independent set of size at least $\alpha(H)/\rho$, that is, a $\rho$-approximation for \textsc{Maximum Independent Set}. Since the latter is NP-hard to approximate within $|W|^{1-\varepsilon}$ for every $\varepsilon>0$~\cite{Hastad1999,Zuckerman2007}, and $|W|=h$, no polynomial-time $h^{1-\varepsilon}$-approximation for Problem~\ref{Pro:1} can exist unless $\mathrm{P}=\mathrm{NP}$.
\end{proof}

Theorem~\ref{thm:inapprox} shows that the linear factor of Theorem~\ref{thm:greedy-ratio} is essentially unavoidable: no efficient algorithm guarantees a ratio better than $h^{1-\varepsilon}$, so the greedy heuristic is, up to the exponent, of the best achievable approximation order. This is what motivates the heuristic and column-generation methods of this section, rather than the pursuit of a constant-factor guarantee.

We note one structural limitation of the algorithm.
The scores $\theta(C)$ used to sort the causalities in Step 1 are computed once on the original graph $G$ and never updated.
As causalities are accepted and more nodes are removed, the original $\theta$ values become stale, and the resulting order of consideration may no longer be the best for the current state.
Recomputing $\theta$ after each acceptance would address this, at the cost of an additional shortest-path computation per accepted causality.
The experiments in Section~\ref{S:6} show that the loss from keeping the static order is small in most cases and concentrated on the harder instances.

\section{Experiments}
\label{S:6}

In this section, we demonstrate the performance of the exact MILP optimization model, the column generation algorithm (CG), and the greedy algorithm (Algorithm 2) in experiments involving real-world networks. All algorithms are implemented in C++ and run on a workstation with an Intel Core i9-14900K CPU and 64 \, GB of RAM. 
We solved the optimization models with both Gurobi~\cite{gurobi} and CPLEX~\cite{cplex} solvers. Since CPLEX was consistently faster in our cases, we report only the CPLEX results.

\subsection{Case Study I: IEEE 300 Network}
\begin{figure*}[!t]
    \centering
    \includegraphics[width=0.6\textwidth]{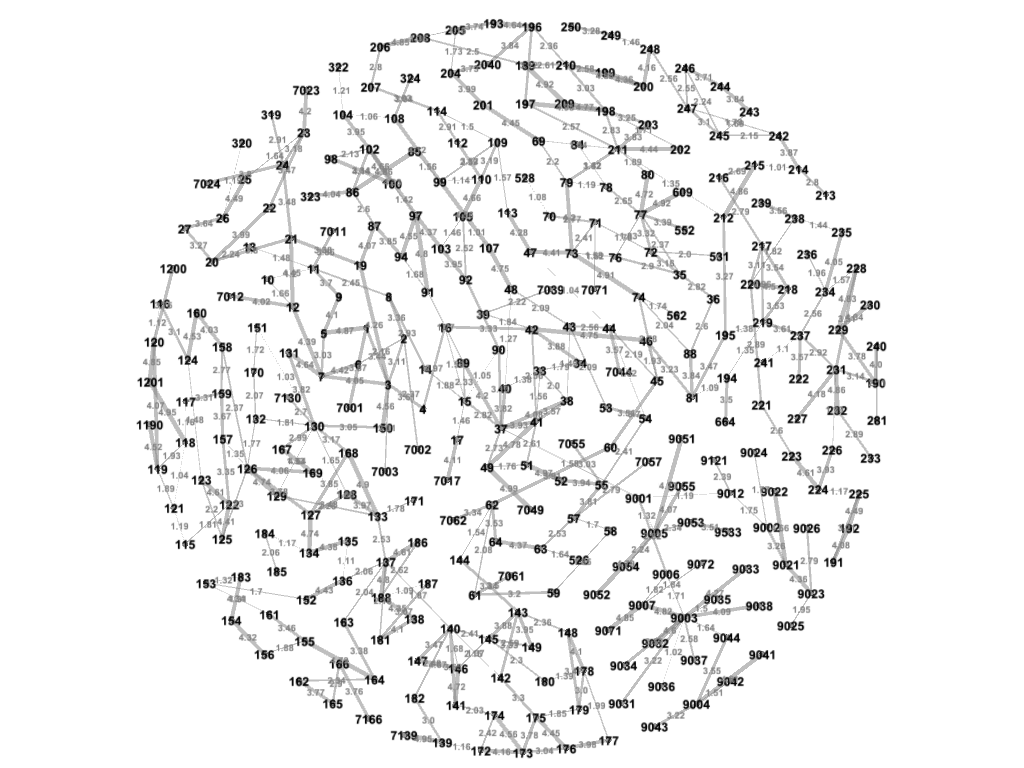}
    \caption{IEEE 300 network with synthetic edge weights.}
    \label{IEEE_300}
\end{figure*}

We explore the connectivity of the IEEE 300 bus system, the topology of which is found in Fig. \ref{IEEE_300}. 
The system has 300 nodes and 409 edges.
We randomly generate 25 causalities as shown in Table \ref{IEEE_300_causality}.
In each causality $C:V_1\Rightarrow V_2$, the number of nodes in $V_1$ or $V_2$ follows a discrete uniform distribution in $\{1, 2, 3, 4, 5\}$.

\begin{table}[!t]
    \centering
    \caption{Causality definition generated randomly for the IEEE 300 network, where numbers in the causality sets represent node labels from the original network data set.}
    \footnotesize
    \begin{tabular}{c|c|c}
        No. & $V_1$ & $V_2$  \\ \hline
         1  &136, 127 &113 \\
         2  &27, 90, 175 &148, 107, 102 \\
         3  &181 &89, 137, 52 \\
         4  &86, 9044 &5, 171 \\
         5  &7130, 42, 9036 &25, 79, 198, 240 \\
         6  &61, 9042, 188 &552, 201, 41 \\
         7  &226, 53, 104, 124 &9044, 241, 228 \\
         8  &243, 88, 203, 248, 176 &526, 44, 102, 129 \\
         9  &196, 117, 48, 39, 235 &145 \\
         10  &9055, 72 &9053, 241, 224, 9026, 201 \\
         11  &121, 126, 161, 178, 9037 &9038, 122, 125, 72 \\
         12  &526, 153, 9055 &16, 225, 5, 109 \\
         13  &108, 72, 179, 1201 &188 \\
         14  &160, 9035, 142, 60, 7139 &140, 9023 \\
         15  &214, 116, 22 &103, 188, 9052, 182, 156 \\
         16  &198, 9071, 237, 23, 166 &9022 \\
         17  &322, 214, 103, 205, 9037 &2040, 9071 \\
         18  &20 &130, 72, 218, 128 \\
         19  &4, 203, 90 &64, 62, 204, 202 \\
         20  &135 &151 \\
         21  &103, 211 &53, 9023, 129, 115, 73 \\
         22  &226, 7055, 169, 121 &141, 229, 227, 120 \\
         23  &25 &9002 \\
         24  &14, 158, 34, 27 &224, 609 \\
         25  &114 &168, 156, 48, 9022, 45 \\ \hline
    \end{tabular}
    \label{IEEE_300_causality}
\end{table}

\subsubsection{Unweighted Network}
We pick four pairs of terminal nodes and compare the exact MILP model, the column generation algorithm (CG), and the greedy algorithm. The results are reported in Table~\ref{result_IEEE_300}, where for each method we list the number of applied causalities $|C|$ and the running time in seconds.

All three methods are fast on this network: the greedy algorithm runs in well under a millisecond, while the MILP and CG both stay below $0.04$ seconds. CG returns exactly the same number of causalities as the MILP on every configuration, confirming that the column generation reaches the optimum on the IEEE 300 instances; their running times are of the same order, with CG occasionally faster (e.g.,\ for $s=2,\,t=7061,\,D=11$, CG takes $0.004$ seconds against $0.022$ for the MILP) and occasionally slower (e.g.,\ $s=3,\,t=7039,\,D=18$, $0.026$ against $0.017$ seconds).

The greedy algorithm matches the optimum in most cases, but not always. For $s=3,\,t=7039$ and $s=7166,\,t=9533$ it finds the optimal count for every value of $D$ ($20$--$22$ and $21$, respectively). The gap appears on the other two pairs: for $s=2,\,t=7061,\,D=13$ greedy applies only $20$ causalities, three fewer than $23$ found by exact methods, and the same shortfall persists at $D=17$ (greedy $21$ against $23$). 
Similarly, for $s=199,\,t=9003$ the greedy stays in the causalities $22$ for all $D$, one to two fewer than the optimum, which itself grows from $23$ to $24$ as $D$ increases from $16$ to $24$.

The optimum is often insensitive to the budget: for $s=7166,\,t=9533$ the number of applied causalities remains $21$ at $D=20,24,28$, because the optimization keeps the shortest $s$--$t$ path unchanged and applies every causality that does not touch a node on that path.

\begin{table*}[!t]
    \centering
    \caption{Experiments for unweighted IEEE 300}
    \label{result_IEEE_300}
    \footnotesize
    \setlength{\tabcolsep}{3pt}
    \begin{tabular}{c|c|c|ccc|ccc|ccc}
    \hline
    \multirow{2}{*}{$\{s,t\}$} & \multirow{2}{*}{$D_0$} & \multirow{2}{*}{$D$} & \multicolumn{3}{c|}{MILP} & \multicolumn{3}{c|}{CG} & \multicolumn{3}{c}{Greedy} \\
     &  &  & $|C|$ & nodes & $t$(s) & $|C|$ & nodes & $t$(s) & $|C|$ & nodes & $t$(s) \\ \hline
    \multirow{3}{*}{$\{3,7039\}$} & \multirow{3}{*}{5} & 6 & 20 & 91 & 0.00567 & 20 & 91 & 0.00288 & 20 & 91 & 0.000441 \\
      &  & 10 & 21 & 93 & 0.0158 & 21 & 93 & 0.0173 & 21 & 93 & 0.000448 \\
      &  & 18 & 22 & 98 & 0.0166 & 22 & 98 & 0.0258 & 22 & 98 & 0.000519 \\
    \hline
    \multirow{3}{*}{$\{2,7061\}$} & \multirow{3}{*}{10} & 11 & 20 & 87 & 0.0218 & 20 & 87 & 0.00361 & 20 & 87 & 0.000821 \\
      &  & 13 & 23 & 103 & 0.0163 & 23 & 103 & 0.00492 & 20 & 87 & 0.000767 \\
      &  & 17 & 23 & 103 & 0.0159 & 23 & 101 & 0.00528 & 21 & 92 & 0.000769 \\
    \hline
    \multirow{3}{*}{$\{199,9003\}$} & \multirow{3}{*}{13} & 16 & 23 & 104 & 0.0175 & 23 & 104 & 0.00466 & 22 & 96 & 0.000641 \\
      &  & 20 & 23 & 104 & 0.0372 & 23 & 104 & 0.00608 & 22 & 96 & 0.000689 \\
      &  & 24 & 24 & 107 & 0.0148 & 24 & 107 & 0.0236 & 22 & 96 & 0.000664 \\
    \hline
    \multirow{3}{*}{$\{7166,9533\}$} & \multirow{3}{*}{19} & 20 & 21 & 93 & 0.0174 & 21 & 93 & 0.00408 & 21 & 93 & 0.000738 \\
      &  & 24 & 21 & 93 & 0.0171 & 21 & 93 & 0.00387 & 21 & 93 & 0.000706 \\
      &  & 28 & 21 & 93 & 0.0183 & 21 & 93 & 0.0203 & 21 & 93 & 0.000714 \\
    \hline
    \end{tabular}
\end{table*}

The number of applied causalities also tells us how many nodes are removed. For $s=199,\,t=9003,\,D=20$, the $23$ applied causalities remove $104$ nodes while keeping the distance between the two terminals at most $20$. This highlights the computational advantage of reasoning in terms of causalities rather than individual nodes: searching directly over node failures would require examining up to $2^{104}$ subsets, whereas the causal model restricts the search to the $25$ causalities of Table~\ref{IEEE_300_causality} (i.e.,\ at most $2^{25}$ combinations).

\subsubsection{Weighted Network}
For the same network we draw a real edge weight uniformly between $1$ and $5$, reuse the four pairs of terminals, and replace the length budget $D$ by three values of the factor $\alpha$. Table~\ref{result_IEEE_300_weighted} reports the outcome.

The picture is consistent with the unweighted case. CG again coincides with the MILP on every instance, and both stay below $0.08$ seconds, so weighting the edges does not appreciably change the difficulty of the IEEE 300 instances. The greedy algorithm is the fastest (a few milliseconds); it is optimal for $s=3,\,t=7039$ and $s=7166,\,t=9533$, but applies two fewer causalities ($21$ instead of $23$) for $s=2,\,t=7061$ at every $\alpha$, and one to two fewer for $s=199,\,t=9003$.

The solutions are largely insensitive to $\alpha$: for $s=7166,\,t=9533$ all three methods return $21$ causalities at $\alpha=1.5,2.5,3.5$, and for $s=2,\,t=7061$ the counts are frozen at $21$ (greedy) and $23$ (exact methods). The exact optimum reacts to $\alpha$ on the other two pairs: it rises from $20$ at $\alpha=1.5$ to $22$ at $\alpha=2.5$ and then plateaus for $s=3,\,t=7039$, and from $23$ to $24$ for $s=199,\,t=9003$.

\begin{table*}[!t]
    \centering
    \caption{Experiments for edge-weighted IEEE 300}
    \label{result_IEEE_300_weighted}
    \footnotesize
    \setlength{\tabcolsep}{3pt}
    \begin{tabular}{c|c|c|ccc|ccc|ccc}
    \hline
    \multirow{2}{*}{$\{s,t\}$} & \multirow{2}{*}{$D_0$} & \multirow{2}{*}{$\alpha$} & \multicolumn{3}{c|}{MILP} & \multicolumn{3}{c|}{CG} & \multicolumn{3}{c}{Greedy} \\
     &  &  & $|C|$ & nodes & $t$(s) & $|C|$ & nodes & $t$(s) & $|C|$ & nodes & $t$(s) \\ \hline
    \multirow{3}{*}{$\{3,7039\}$} & \multirow{3}{*}{15.10} & 1.5 & 20 & 91 & 0.0296 & 20 & 91 & 0.00453 & 20 & 91 & 0.00152 \\
      &  & 2.5 & 22 & 98 & 0.0169 & 22 & 98 & 0.027 & 22 & 98 & 0.00148 \\
      &  & 3.5 & 22 & 98 & 0.0143 & 22 & 98 & 0.0348 & 22 & 98 & 0.00131 \\
    \hline
    \multirow{3}{*}{$\{2,7061\}$} & \multirow{3}{*}{32.63} & 1.5 & 23 & 103 & 0.0163 & 23 & 103 & 0.00535 & 21 & 92 & 0.00287 \\
      &  & 2.5 & 23 & 103 & 0.0156 & 23 & 103 & 0.00817 & 21 & 92 & 0.00262 \\
      &  & 4.0 & 23 & 101 & 0.0155 & 23 & 103 & 0.00779 & 21 & 92 & 0.00269 \\
    \hline
    \multirow{3}{*}{$\{199,9003\}$} & \multirow{3}{*}{31.52} & 1.5 & 23 & 103 & 0.0361 & 23 & 103 & 0.00928 & 22 & 96 & 0.00162 \\
      &  & 2.5 & 24 & 107 & 0.017 & 24 & 107 & 0.0597 & 22 & 96 & 0.0016 \\
      &  & 4.0 & 24 & 107 & 0.0175 & 24 & 107 & 0.0783 & 22 & 96 & 0.0017 \\
    \hline
    \multirow{3}{*}{$\{7166,9533\}$} & \multirow{3}{*}{54.98} & 1.5 & 21 & 93 & 0.0181 & 21 & 93 & 0.00553 & 21 & 93 & 0.00233 \\
      &  & 2.5 & 21 & 93 & 0.0176 & 21 & 93 & 0.00863 & 21 & 93 & 0.00219 \\
      &  & 3.5 & 21 & 93 & 0.0197 & 21 & 93 & 0.00911 & 21 & 93 & 0.00242 \\
    \hline
    \end{tabular}
\end{table*}

\subsection{Case Study II: Power 5000}

A power network with 4945 nodes and 6570 edges is generated using the method in \cite{Soltan} and introduced to investigate the scalability of the proposed algorithms.
We also randomly generate 40 causalities.
For each causality, $C: V_1 \Rightarrow V_2$, the number of nodes in $V_1$ or $V_2$ is sampled from a discrete uniform distribution between 1 and 60.

\subsubsection{Unweighted Network}
Each edge has unit weight. We select four pairs of terminals and report in Table~\ref{result_unweighted_power_5000} the number of applied causalities and the running time of the three methods.

On this larger network the value of CG becomes evident. CG returns exactly the same number of causalities as the exact MILP on all twenty configurations, yet it is roughly one order of magnitude faster: the MILP spends between $0.6$ and $1.1$ seconds per instance, whereas CG ranges from $0.02$ to $0.44$ seconds (for example, for $s=184,\,t=187,\,D=8$ the MILP takes $0.71$ seconds against $0.025$ for CG). The greedy algorithm remains the cheapest, always under $0.04$ seconds.

In terms of quality, the greedy heuristic is again close to the optimum but loses a few causalities on the harder pairs: for $s=709,\,t=1314,\,D=28$ it applies $27$ causalities against $29$ for the exact methods, and for $s=1035,\,t=1200$ it trails by one unit at $D=28,30,32$. The most striking instance is $s=709,\,t=1314$ at $D=20$ and $D=22$: \emph{no} causality can be applied by any method ($|C|=0$), because every available causality, once its closure is removed, pushes the distance between the two terminals beyond the budget; only from $D=24$ onwards does the network admit causalities, jumping abruptly to $26$ and then to $29$. At the opposite end, the pair $s=10,\,t=14$ is completely insensitive to the budget, applying $36$ causalities for every $D$ from $14$ to $22$.

The number of applied causalities also measures how many nodes are removed: for $s=184,\,t=187$ the $36$ causalities knock out $1626$ nodes, roughly a third of the $4945$-node network, while keeping the two terminals within the prescribed distance.

\begin{table*}[!t]
    \centering
    \caption{Experiments for unweighted Power 5000}
    \label{result_unweighted_power_5000}
    \footnotesize
    \setlength{\tabcolsep}{3pt}
    \begin{tabular}{c|c|c|ccc|ccc|ccc}
    \hline
    \multirow{2}{*}{$\{s,t\}$} & \multirow{2}{*}{$D_0$} & \multirow{2}{*}{$D$} & \multicolumn{3}{c|}{MILP} & \multicolumn{3}{c|}{CG} & \multicolumn{3}{c}{Greedy} \\
     &  &  & $|C|$ & nodes & $t$(s) & $|C|$ & nodes & $t$(s) & $|C|$ & nodes & $t$(s) \\ \hline
    \multirow{5}{*}{$\{184,187\}$} & \multirow{5}{*}{6} & 8 & 36 & 1626 & 0.712 & 36 & 1626 & 0.0254 & 36 & 1626 & 0.00461 \\
      &  & 10 & 36 & 1626 & 0.639 & 36 & 1626 & 0.0204 & 36 & 1626 & 0.0035 \\
      &  & 12 & 36 & 1626 & 0.666 & 36 & 1626 & 0.0232 & 36 & 1626 & 0.0034 \\
      &  & 14 & 37 & 1695 & 0.778 & 37 & 1695 & 0.0305 & 36 & 1626 & 0.00341 \\
      &  & 16 & 38 & 1728 & 0.607 & 38 & 1728 & 0.443 & 38 & 1728 & 0.00563 \\
    \hline
    \multirow{5}{*}{$\{10,14\}$} & \multirow{5}{*}{9} & 14 & 36 & 1602 & 0.706 & 36 & 1602 & 0.0383 & 36 & 1602 & 0.00466 \\
      &  & 16 & 36 & 1602 & 0.682 & 36 & 1602 & 0.0511 & 36 & 1602 & 0.00382 \\
      &  & 18 & 36 & 1595 & 0.793 & 36 & 1602 & 0.0665 & 36 & 1602 & 0.0037 \\
      &  & 20 & 36 & 1595 & 0.667 & 36 & 1595 & 0.115 & 36 & 1602 & 0.0043 \\
      &  & 22 & 36 & 1595 & 0.669 & 36 & 1595 & 0.102 & 36 & 1602 & 0.00385 \\
    \hline
    \multirow{5}{*}{$\{709,1314\}$} & \multirow{5}{*}{24} & 20 & 0 & 0 & 0.944 & 0 & 0 & 0.000647 & 0 & 0 & 0.0236 \\
      &  & 22 & 0 & 0 & 0.881 & 0 & 0 & 0.000676 & 0 & 0 & 0.0224 \\
      &  & 24 & 26 & 1204 & 0.86 & 26 & 1204 & 0.0189 & 26 & 1204 & 0.0184 \\
      &  & 26 & 28 & 1317 & 0.947 & 28 & 1317 & 0.0241 & 27 & 1258 & 0.0202 \\
      &  & 28 & 29 & 1372 & 1.09 & 29 & 1372 & 0.0347 & 27 & 1258 & 0.0183 \\
    \hline
    \multirow{5}{*}{$\{1035,1200\}$} & \multirow{5}{*}{25} & 26 & 31 & 1395 & 0.855 & 31 & 1395 & 0.0245 & 31 & 1395 & 0.0281 \\
      &  & 28 & 33 & 1508 & 0.781 & 33 & 1508 & 0.026 & 32 & 1452 & 0.0359 \\
      &  & 30 & 33 & 1508 & 0.835 & 33 & 1508 & 0.0413 & 32 & 1452 & 0.0275 \\
      &  & 32 & 34 & 1565 & 0.951 & 34 & 1547 & 0.0671 & 33 & 1497 & 0.0312 \\
      &  & 34 & 35 & 1610 & 1.12 & 35 & 1610 & 0.414 & 35 & 1610 & 0.0272 \\
    \hline
    \end{tabular}
\end{table*}

\subsubsection{Weighted Network}
We assign random weights between $1$ and $5$ to the edges and repeat the experiment with several values of $\alpha$ per pair. Table~\ref{result_weighted_power_5000} collects the results.

The running times grow with respect to the unweighted case --- the MILP now reaches $6.5$ seconds for $s=1035,\,t=1200$ at $\alpha=2$ --- and CG keeps a clear advantage, solving the same instance in about $1.1$ seconds. More interestingly, this is the only setting in which CG no longer always matches the exact optimum. CG stays optimal for $s=184,\,t=187$ (all configurations) and for $s=10,\,t=14$, where it reaches $36$ causalities against the greedy's $33$--$34$; but on the two hardest pairs it occasionally applies fewer causalities than the MILP: for $s=709,\,t=1314$ it drops to $30$ at $\alpha=3,3.5,4$ where the MILP applies $33$, and for $s=1035,\,t=1200$ at $\alpha=2.5$ it returns $32$ against the MILP's $35$. This is the expected consequence of solving the final integer master only over the pool of generated columns: on these instances the pool does not contain the paths needed to certify the integer optimum. In a few of these cases even the greedy heuristic beats CG (for example $s=709,\,t=1314,\,\alpha=3$: greedy $32$, MILP $33$, CG $30$), which suggests enlarging the column pool or embedding the pricing in a branch-and-price scheme on such instances.

The weighted version of $s=709,\,t=1314$ is particularly telling: whereas in the unweighted case no causality could be applied below $D=24$, weighting the edges produces several $s$--$t$ paths of similar length, so the network can absorb more causalities before the distance constraint is violated, and the applied count rises from $0$ (unweighted, low budget) to around $30$--$33$.

\begin{table*}[!t]
    \centering
    \caption{Experiments for edge-weighted Power 5000}
    \label{result_weighted_power_5000}
    \footnotesize
    \setlength{\tabcolsep}{3pt}
    \begin{tabular}{c|c|c|ccc|ccc|ccc}
    \hline
    \multirow{2}{*}{$\{s,t\}$} & \multirow{2}{*}{$D_0$} & \multirow{2}{*}{$\alpha$} & \multicolumn{3}{c|}{MILP} & \multicolumn{3}{c|}{CG} & \multicolumn{3}{c}{Greedy} \\
     &  &  & $|C|$ & nodes & $t$(s) & $|C|$ & nodes & $t$(s) & $|C|$ & nodes & $t$(s) \\ \hline
    \multirow{6}{*}{$\{184,187\}$} & \multirow{6}{*}{23.58} & 1.5 & 36 & 1626 & 0.91 & 36 & 1626 & 0.0214 & 36 & 1626 & 0.0138 \\
      &  & 2.0 & 36 & 1626 & 0.951 & 36 & 1626 & 0.0239 & 36 & 1626 & 0.0163 \\
      &  & 2.5 & 36 & 1587 & 0.887 & 36 & 1626 & 0.0264 & 36 & 1626 & 0.0132 \\
      &  & 3.0 & 38 & 1728 & 0.807 & 38 & 1728 & 0.368 & 38 & 1728 & 0.013 \\
      &  & 3.5 & 38 & 1707 & 0.577 & 38 & 1728 & 0.378 & 38 & 1728 & 0.0131 \\
      &  & 4.0 & 38 & 1707 & 0.591 & 38 & 1728 & 0.387 & 38 & 1728 & 0.013 \\
    \hline
    \multirow{6}{*}{$\{10,14\}$} & \multirow{6}{*}{43.74} & 1.5 & 36 & 1602 & 0.7 & 36 & 1602 & 0.0271 & 33 & 1439 & 0.0237 \\
      &  & 2.0 & 36 & 1602 & 0.768 & 36 & 1602 & 0.058 & 33 & 1439 & 0.0238 \\
      &  & 2.5 & 36 & 1595 & 0.68 & 36 & 1602 & 0.0921 & 34 & 1485 & 0.0236 \\
      &  & 3.0 & 36 & 1595 & 0.626 & 36 & 1602 & 0.0918 & 34 & 1485 & 0.0235 \\
      &  & 3.5 & 36 & 1595 & 0.625 & 36 & 1602 & 0.101 & 34 & 1485 & 0.0239 \\
      &  & 4.0 & 36 & 1595 & 0.624 & 36 & 1602 & 0.106 & 34 & 1485 & 0.024 \\
    \hline
    \multirow{6}{*}{$\{709,1314\}$} & \multirow{6}{*}{95.56} & 1.5 & 31 & 1399 & 0.915 & 31 & 1399 & 0.532 & 30 & 1368 & 0.0747 \\
      &  & 2.0 & 32 & 1453 & 1.43 & 31 & 1365 & 0.265 & 32 & 1442 & 0.0716 \\
      &  & 2.5 & 32 & 1453 & 1.67 & 32 & 1398 & 1.98 & 32 & 1442 & 0.0742 \\
      &  & 3.0 & 33 & 1439 & 2.56 & 30 & 1334 & 1.45 & 32 & 1442 & 0.0791 \\
      &  & 3.5 & 33 & 1439 & 2.14 & 30 & 1276 & 1.04 & 32 & 1442 & 0.0708 \\
      &  & 4.0 & 33 & 1439 & 2.2 & 30 & 1306 & 1.91 & 32 & 1442 & 0.0792 \\
    \hline
    \multirow{6}{*}{$\{1035,1200\}$} & \multirow{6}{*}{123.47} & 1.5 & 35 & 1610 & 0.999 & 35 & 1610 & 0.635 & 35 & 1610 & 0.119 \\
      &  & 2.0 & 35 & 1583 & 6.5 & 33 & 1464 & 1.09 & 34 & 1529 & 0.12 \\
      &  & 2.5 & 35 & 1583 & 5.83 & 32 & 1429 & 1.57 & 35 & 1568 & 0.119 \\
      &  & 3.0 & 36 & 1614 & 5.59 & 33 & 1489 & 2.34 & 35 & 1568 & 0.122 \\
      &  & 3.5 & 36 & 1595 & 1.56 & 34 & 1600 & 2 & 35 & 1568 & 0.125 \\
      &  & 4.0 & 36 & 1607 & 3.01 & 34 & 1600 & 2.12 & 35 & 1568 & 0.119 \\
    \hline
    \end{tabular}
\end{table*}

\subsubsection{Scalability Analysis}
Comparing the IEEE 300 and the Power 5000 networks shows how the three methods scale. The exact MILP is the most affected: its average solving time grows from under $0.02$ seconds on IEEE 300 to about one second --- and up to $6.5$ seconds in the weighted case --- on Power 5000. The greedy algorithm scales gracefully, rising only from a fraction of a millisecond to at most about a tenth of a second. CG lies in between, but much closer to the greedy than to the MILP: it stays between a few hundredths of a second and roughly two seconds on Power 5000 (i.e.,\ a small fraction of the MILP time), while returning the same solution (i.e., on a few weighted instances, a solution that is only one to three causalities short). CG therefore offers the most favorable accuracy/time trade-off as the network grows, and the greedy algorithm remains a strong alternative whenever a small loss of optimality is tolerable.

\subsubsection{Sensitivity Analysis}
Across both networks the solutions are remarkably stable with respect to the budget. Once the source and the terminal are fixed, the number of applied causalities is essentially determined and changes little as $D$ or $\alpha$ varies: see $s=7166,\,t=9533$ on IEEE 300 (constant at $21$) and $s=10,\,t=14$ on Power 5000 (constant at $36$). The reason is that, for a small budget, the network already applies every causality that does not lengthen the $s$--$t$ path beyond the limit; enlarging the budget rarely unlocks new causalities, because the remaining ones would disconnect $s$ from $t$ or push their distance past the bound. The clearest exception is the threshold behaviour of $s=709,\,t=1314$ on the unweighted Power 5000, where the count jumps from $0$ to $26$ as soon as $D$ reaches $24$.

\subsection{Case Study III: Large-Scale Synthetic Networks}
To stress-test the three methods on much larger instances, we generate synthetic Erd\H{o}s--R\'enyi random graphs~\cite{ErdosRenyi1960,Gilbert1959} and keep their largest connected component: \emph{sparse} graphs with edge probability $p=0.001$ (average degree $\approx 10$) on $10{,}000$ and $15{,}000$ nodes, and a \emph{dense} graph with $p=0.01$ (average degree $\approx 100$) on $10{,}000$ nodes. On each graph we generate causalities with two settings of increasing difficulty: a base set with up to $120$ nodes per side and $160$ causalities ($\text{bound}=120,\ h=160$), and a harder set with $\text{bound}=180,\ h=240$. These instances are deliberately large and rich in causalities, so as to drive the exact MILP to its time limit and thereby expose the need for scalable alternatives.

On the larger and denser of these instances the MILP can no longer certify optimality within its time budget. Wherever the MILP does solve to optimality, we use its certified optimum as the reference; on the harder instances, where no method certifies optimality, we compare the three methods directly by the number of applicable causalities they find and by their running time. Since the objective is to maximise the number of applied causalities, the method returning the largest feasible count is the best solution found.

\subsubsection{Sparse 10k Networks}
Tables~\ref{tab:syn_sparse10k_u} and~\ref{tab:syn_sparse10k_w} report the sparse $10{,}000$-node instances. Here the MILP still solves every instance to optimality, so we compare the three methods against its certified optimum. On the unweighted graph (Table~\ref{tab:syn_sparse10k_u}) the budget is essentially non-binding at these large values of $D$, so the number of applied causalities saturates around $152$--$156$ of the $160$ available and barely moves as $D$ grows from $15$ to $5000$. The three methods agree to within at most one or two causalities on every configuration; what separates them is speed: the MILP is slow and erratic, ranging from about $5$ seconds for $s=1000,\,t=8000$ up to a $275$-second spike for $s=500,\,t=5000$ at $D=15$, whereas column generation stays under $8$ seconds and the greedy heuristic under one second.

\begin{table}[!t]
    \centering
    \caption{Synthetic sparse 10k network, unweighted (bound $=120$, $h=160$)}
    \label{tab:syn_sparse10k_u}
    \footnotesize
    \setlength{\tabcolsep}{3pt}
    \begin{tabular}{c|c|c|cc|cc|cc}
    \hline
    \multirow{2}{*}{$\{s,t\}$} & \multirow{2}{*}{$D_0$} & \multirow{2}{*}{$D$} & \multicolumn{2}{c|}{MILP} & \multicolumn{2}{c|}{CG} & \multicolumn{2}{c}{Greedy} \\
     &  &  & $|C|$ & $t$(s) & $|C|$ & $t$(s) & $|C|$ & $t$(s) \\ \hline
    \multirow{5}{*}{$\{100,9900\}$} & \multirow{5}{*}{5} & 15 & 154 & 9.39 & 154 & 3.9 & 154 & 0.748 \\
      &  & 50 & 154 & 10.1 & 154 & 3.28 & 154 & 0.76 \\
      &  & 500 & 154 & 10 & 154 & 3.25 & 154 & 0.751 \\
      &  & 1000 & 154 & 10.4 & 154 & 3.25 & 154 & 0.766 \\
      &  & 5000 & 154 & 10.6 & 154 & 3.28 & 154 & 0.753 \\
    \hline
    \multirow{5}{*}{$\{500,5000\}$} & \multirow{5}{*}{4} & 15 & 154 & 275 & 154 & 7.55 & 153 & 0.296 \\
      &  & 50 & 155 & 57.1 & 153 & 7.34 & 154 & 0.3 \\
      &  & 500 & 155 & 56.8 & 153 & 7.4 & 154 & 0.31 \\
      &  & 1000 & 155 & 58.2 & 153 & 7.26 & 154 & 0.296 \\
      &  & 5000 & 155 & 57.6 & 153 & 7.33 & 154 & 0.305 \\
    \hline
    \multirow{5}{*}{$\{1000,8000\}$} & \multirow{5}{*}{5} & 15 & 153 & 5.99 & 153 & 1.41 & 152 & 0.89 \\
      &  & 50 & 153 & 5.83 & 153 & 1.39 & 153 & 0.883 \\
      &  & 500 & 153 & 5.03 & 153 & 1.39 & 153 & 0.862 \\
      &  & 1000 & 153 & 5.1 & 153 & 1.42 & 153 & 0.871 \\
      &  & 5000 & 153 & 4.99 & 153 & 1.42 & 153 & 0.883 \\
    \hline
    \multirow{5}{*}{$\{2500,7500\}$} & \multirow{5}{*}{5} & 15 & 156 & 46.2 & 155 & 2.45 & 155 & 0.725 \\
      &  & 50 & 156 & 40.7 & 155 & 2.5 & 156 & 0.748 \\
      &  & 500 & 156 & 39.5 & 155 & 2.39 & 156 & 0.727 \\
      &  & 1000 & 156 & 39.2 & 155 & 2.39 & 156 & 0.757 \\
      &  & 5000 & 156 & 39.3 & 155 & 2.41 & 156 & 0.733 \\
    \hline
    \end{tabular}
\end{table}

Weighting the edges (Table~\ref{tab:syn_sparse10k_w}) raises every running time without changing the picture. The MILP now needs up to about $55$ seconds (for $s=2500,\,t=7500$), column generation stays around $9$ seconds and the greedy under $6$, and the applied counts are frozen across all values of $\alpha$ --- $154$ for $s=100,\,t=9900$, $153$ for $s=1000,\,t=8000$ --- confirming that the length budget is slack. The three methods again return almost identical counts, differing by at most two causalities, with column generation and the greedy an order of magnitude faster than the exact model.

\begin{table}[!t]
    \centering
    \caption{Synthetic sparse 10k network, edge-weighted (bound $=120$, $h=160$)}
    \label{tab:syn_sparse10k_w}
    \footnotesize
    \setlength{\tabcolsep}{3pt}
    \begin{tabular}{c|c|c|cc|cc|cc}
    \hline
    \multirow{2}{*}{$\{s,t\}$} & \multirow{2}{*}{$D_0$} & \multirow{2}{*}{$\alpha$} & \multicolumn{2}{c|}{MILP} & \multicolumn{2}{c|}{CG} & \multicolumn{2}{c}{Greedy} \\
     &  &  & $|C|$ & $t$(s) & $|C|$ & $t$(s) & $|C|$ & $t$(s) \\ \hline
    \multirow{4}{*}{$\{100,9900\}$} & \multirow{4}{*}{9.98} & 15 & 154 & 45.2 & 154 & 8.89 & 154 & 2.3 \\
      &  & 50 & 154 & 48.3 & 154 & 9 & 154 & 2.31 \\
      &  & 500 & 154 & 45.9 & 154 & 8.98 & 154 & 2.3 \\
      &  & 800 & 154 & 46.3 & 154 & 9.01 & 154 & 2.31 \\
    \hline
    \multirow{4}{*}{$\{500,5000\}$} & \multirow{4}{*}{9.18} & 15 & 155 & 23.9 & 153 & 8.96 & 154 & 0.737 \\
      &  & 50 & 155 & 24.8 & 153 & 9.07 & 154 & 0.748 \\
      &  & 500 & 155 & 25 & 153 & 8.99 & 154 & 0.747 \\
      &  & 800 & 155 & 24.9 & 153 & 8.89 & 154 & 0.73 \\
    \hline
    \multirow{4}{*}{$\{1000,8000\}$} & \multirow{4}{*}{11.50} & 15 & 153 & 6.55 & 153 & 0.36 & 153 & 3.98 \\
      &  & 50 & 153 & 6.51 & 153 & 0.357 & 153 & 3.99 \\
      &  & 500 & 153 & 8.76 & 153 & 0.409 & 153 & 3.97 \\
      &  & 800 & 153 & 8.75 & 153 & 0.477 & 153 & 5.5 \\
    \hline
    \multirow{4}{*}{$\{2500,7500\}$} & \multirow{4}{*}{11.96} & 15 & 156 & 54.7 & 155 & 4.7 & 154 & 5.19 \\
      &  & 50 & 156 & 51.5 & 155 & 4.44 & 154 & 5.39 \\
      &  & 500 & 156 & 51.1 & 155 & 4.94 & 154 & 4.95 \\
      &  & 800 & 156 & 49.5 & 155 & 3.55 & 154 & 5.54 \\
    \hline
    \end{tabular}
\end{table}

\subsubsection{Dense 10k Networks}
The dense graph is where the exact model finally breaks down (Tables~\ref{tab:syn_dense10k_u_b120}--\ref{tab:syn_dense10k_w_b180}). We discuss the four configurations in turn, from the base causality set ($\text{bound}=120$) to the harder one ($\text{bound}=180$).

On the unweighted base instances (Table~\ref{tab:syn_dense10k_u_b120}), the MILP still reaches the optimum on every pair, certifying $154$ to $157$ causalities, exactly as the greedy heuristic and column generation. The only penalty is speed: the MILP fluctuates between about $20$ and $79$ seconds, while column generation is uniformly around $3.5$ seconds and the greedy around $5$--$7$. The budget is again slack, so the counts are constant across $D$ for every pair. This is the last dense setting in which the exact model is fully reliable.

\begin{table}[!t]
    \centering
    \caption{Synthetic dense 10k network, unweighted (bound $=120$, $h=160$)}
    \label{tab:syn_dense10k_u_b120}
    \footnotesize
    \setlength{\tabcolsep}{3pt}
    \begin{tabular}{c|c|c|cc|cc|cc}
    \hline
    \multirow{2}{*}{$\{s,t\}$} & \multirow{2}{*}{$D_0$} & \multirow{2}{*}{$D$} & \multicolumn{2}{c|}{MILP} & \multicolumn{2}{c|}{CG} & \multicolumn{2}{c}{Greedy} \\
     &  &  & $|C|$ & $t$(s) & $|C|$ & $t$(s) & $|C|$ & $t$(s) \\ \hline
    \multirow{5}{*}{$\{100,9900\}$} & \multirow{5}{*}{3} & 15 & 155 & 60.3 & 155 & 3.68 & 155 & 7.31 \\
      &  & 50 & 155 & 21.8 & 155 & 3.72 & 155 & 7.01 \\
      &  & 500 & 155 & 47.8 & 155 & 3.56 & 155 & 7.04 \\
      &  & 1000 & 155 & 47.1 & 155 & 3.72 & 155 & 7.18 \\
      &  & 5000 & 155 & 47.4 & 155 & 3.67 & 155 & 6.97 \\
    \hline
    \multirow{5}{*}{$\{500,5000\}$} & \multirow{5}{*}{2} & 15 & 156 & 20.8 & 156 & 3.58 & 156 & 5.84 \\
      &  & 50 & 156 & 20.8 & 156 & 3.4 & 156 & 5.74 \\
      &  & 500 & 156 & 22.3 & 156 & 3.43 & 156 & 5.64 \\
      &  & 1000 & 156 & 74.3 & 156 & 3.43 & 156 & 5.66 \\
      &  & 5000 & 156 & 73.5 & 156 & 3.43 & 156 & 5.68 \\
    \hline
    \multirow{5}{*}{$\{1000,8000\}$} & \multirow{5}{*}{2} & 15 & 154 & 20.7 & 154 & 3.62 & 154 & 4.96 \\
      &  & 50 & 154 & 20.4 & 154 & 3.58 & 154 & 4.97 \\
      &  & 500 & 154 & 55.4 & 154 & 3.58 & 154 & 4.95 \\
      &  & 1000 & 154 & 62.9 & 154 & 3.63 & 154 & 4.94 \\
      &  & 5000 & 154 & 55.8 & 154 & 3.68 & 154 & 5.01 \\
    \hline
    \multirow{5}{*}{$\{2500,7500\}$} & \multirow{5}{*}{2} & 15 & 157 & 57.9 & 157 & 3.52 & 157 & 6.62 \\
      &  & 50 & 157 & 57.2 & 157 & 3.46 & 157 & 6.57 \\
      &  & 500 & 157 & 79 & 157 & 3.57 & 157 & 6.71 \\
      &  & 1000 & 157 & 70.2 & 157 & 3.53 & 157 & 6.64 \\
      &  & 5000 & 157 & 70.2 & 157 & 3.43 & 157 & 6.65 \\
    \hline
    \end{tabular}
\end{table}

Weighting the edges (Table~\ref{tab:syn_dense10k_w_b120}) leaves the greedy and column generation stable at the highest counts on all 20 configurations, but the MILP begins to fail. It holds for $s=1000,\,t=8000$ ($154$ throughout), yet it degrades as $\alpha$ grows on the other pairs, dropping to $136$--$139$ for $s=100,\,t=9900$ and to $142$ for $s=500,\,t=5000$, and it collapses for $s=2500,\,t=7500$, where it returns only $20$, $22$, and $27$ causalities against the $157$ found by both heuristics --- an error of nearly a factor of eight. The MILP now spends its full budget of about $122$ seconds on these pairs, whereas column generation keeps to $5$--$7$ seconds; the greedy is correct but no longer cheap, costing up to $43$ seconds.

\begin{table}[!t]
    \centering
    \caption{Synthetic dense 10k network, edge-weighted (bound $=120$, $h=160$)}
    \label{tab:syn_dense10k_w_b120}
    \footnotesize
    \setlength{\tabcolsep}{3pt}
    \begin{tabular}{c|c|c|cc|cc|cc}
    \hline
    \multirow{2}{*}{$\{s,t\}$} & \multirow{2}{*}{$D_0$} & \multirow{2}{*}{$\alpha$} & \multicolumn{2}{c|}{MILP} & \multicolumn{2}{c|}{CG} & \multicolumn{2}{c}{Greedy} \\
     &  &  & $|C|$ & $t$(s) & $|C|$ & $t$(s) & $|C|$ & $t$(s) \\ \hline
    \multirow{5}{*}{$\{100,9900\}$} & \multirow{5}{*}{4.43} & 15 & 155 & 117 & 155 & 7.24 & 155 & 42.8 \\
      &  & 50 & 139 & 122 & 155 & 7.2 & 155 & 42.5 \\
      &  & 500 & 155 & 121 & 155 & 6.06 & 155 & 42.5 \\
      &  & 1000 & 136 & 123 & 155 & 6 & 155 & 33.1 \\
      &  & 5000 & 136 & 123 & 155 & 6.1 & 155 & 33.7 \\
    \hline
    \multirow{5}{*}{$\{500,5000\}$} & \multirow{5}{*}{3.30} & 15 & 156 & 95 & 156 & 5.04 & 156 & 4.31 \\
      &  & 50 & 156 & 94.6 & 156 & 4.98 & 156 & 4.27 \\
      &  & 500 & 154 & 123 & 156 & 5.17 & 156 & 4.25 \\
      &  & 1000 & 142 & 123 & 156 & 5.1 & 156 & 4.57 \\
      &  & 5000 & 142 & 123 & 156 & 5.17 & 156 & 4.31 \\
    \hline
    \multirow{5}{*}{$\{1000,8000\}$} & \multirow{5}{*}{4.06} & 15 & 154 & 94.3 & 154 & 5.78 & 154 & 14.4 \\
      &  & 50 & 154 & 92.2 & 154 & 5.78 & 154 & 14.3 \\
      &  & 500 & 154 & 92.8 & 154 & 5.79 & 154 & 14.2 \\
      &  & 1000 & 154 & 92.7 & 154 & 5.7 & 154 & 14.2 \\
      &  & 5000 & 154 & 92.8 & 154 & 5.73 & 154 & 14.1 \\
    \hline
    \multirow{5}{*}{$\{2500,7500\}$} & \multirow{5}{*}{3.84} & 15 & 157 & 97.6 & 157 & 4.87 & 157 & 17.6 \\
      &  & 50 & 127 & 122 & 157 & 4.97 & 157 & 17.6 \\
      &  & 500 & 20 & 122 & 157 & 6.34 & 157 & 17.8 \\
      &  & 1000 & 22 & 124 & 157 & 6.23 & 157 & 23.3 \\
      &  & 5000 & 27 & 123 & 157 & 6.31 & 157 & 21.9 \\
    \hline
    \end{tabular}
\end{table}

Enlarging the causality set to $\text{bound}=180$ (Table~\ref{tab:syn_dense10k_u_b180}) pushes the MILP to its five-minute limit (about $303$ seconds) on most pairs and makes its output erratic. It stays within a few causalities of the greedy heuristic and column generation on the easier $s=100,\,t=9900$ and $s=1000,\,t=8000$, but collapses to $105$ and then $22$ causalities for $s=2500,\,t=7500$, against the $232$ returned by both the greedy heuristic and column generation. For $s=500,\,t=5000$ the greedy and column generation stay around $224$--$226$ throughout, while the MILP loses ground once it hits the time limit, dropping to $223$. Column generation stays within one to two causalities of the best result on every pair, in $11$ to $67$ seconds.

\begin{table}[!t]
    \centering
    \caption{Synthetic dense 10k network, unweighted (bound $=180$, $h=240$)}
    \label{tab:syn_dense10k_u_b180}
    \footnotesize
    \setlength{\tabcolsep}{3pt}
    \begin{tabular}{c|c|c|cc|cc|cc}
    \hline
    \multirow{2}{*}{$\{s,t\}$} & \multirow{2}{*}{$D_0$} & \multirow{2}{*}{$D$} & \multicolumn{2}{c|}{MILP} & \multicolumn{2}{c|}{CG} & \multicolumn{2}{c}{Greedy} \\
     &  &  & $|C|$ & $t$(s) & $|C|$ & $t$(s) & $|C|$ & $t$(s) \\ \hline
    \multirow{5}{*}{$\{100,9900\}$} & \multirow{5}{*}{3} & 15 & 226 & 303 & 225 & 63.2 & 225 & 10.8 \\
      &  & 50 & 225 & 303 & 225 & 64.2 & 224 & 14 \\
      &  & 500 & 224 & 302 & 225 & 64.1 & 224 & 14.2 \\
      &  & 1000 & 224 & 303 & 225 & 64.1 & 224 & 14.1 \\
      &  & 5000 & 224 & 304 & 225 & 64 & 224 & 13.9 \\
    \hline
    \multirow{5}{*}{$\{500,5000\}$} & \multirow{5}{*}{2} & 15 & 226 & 25.6 & 224 & 61.5 & 226 & 8.86 \\
      &  & 50 & 226 & 26.8 & 224 & 62 & 226 & 8.81 \\
      &  & 500 & 223 & 303 & 224 & 61.7 & 226 & 8.76 \\
      &  & 1000 & 223 & 303 & 224 & 61.9 & 226 & 8.87 \\
      &  & 5000 & 223 & 304 & 223 & 62.8 & 226 & 11.8 \\
    \hline
    \multirow{5}{*}{$\{1000,8000\}$} & \multirow{5}{*}{2} & 15 & 224 & 302 & 226 & 66.9 & 225 & 10.4 \\
      &  & 50 & 226 & 303 & 226 & 63.7 & 225 & 10.2 \\
      &  & 500 & 225 & 303 & 226 & 62.6 & 225 & 10.1 \\
      &  & 1000 & 225 & 303 & 226 & 66.2 & 225 & 10.5 \\
      &  & 5000 & 225 & 304 & 226 & 65.9 & 225 & 10.4 \\
    \hline
    \multirow{5}{*}{$\{2500,7500\}$} & \multirow{5}{*}{2} & 15 & 232 & 304 & 232 & 10.7 & 232 & 13.3 \\
      &  & 50 & 232 & 303 & 232 & 11.2 & 232 & 13.8 \\
      &  & 500 & 105 & 305 & 232 & 11.4 & 232 & 13.4 \\
      &  & 1000 & 22 & 304 & 232 & 11.6 & 232 & 13.6 \\
      &  & 5000 & 22 & 304 & 232 & 12 & 232 & 14.1 \\
    \hline
    \end{tabular}
\end{table}

The weighted harder set (Table~\ref{tab:syn_dense10k_w_b180}) is the most adversarial. The MILP collapses on three of the four pairs, falling to $12$--$25$ for $s=1000,\,t=8000$ and to $17$--$22$ for $s=2500,\,t=7500$, against the $226$ and $232$ of column generation. Even the greedy heuristic breaks here: for $s=100,\,t=9900$ it exhausts its own time budget ($60.5$ seconds) and returns $0$ causalities, so among the heuristics only column generation remains usable on that pair ($224$); the MILP matches it ($224$) at $\alpha=15$ but degrades to $129$--$138$ at larger $\alpha$. Throughout, column generation is the only method that never collapses, returning between $222$ and $232$ causalities on every pair in $41$ to $65$ seconds.

\begin{table}[!t]
    \centering
    \caption{Synthetic dense 10k network, edge-weighted (bound $=180$, $h=240$)}
    \label{tab:syn_dense10k_w_b180}
    \footnotesize
    \setlength{\tabcolsep}{3pt}
    \begin{tabular}{c|c|c|cc|cc|cc}
    \hline
    \multirow{2}{*}{$\{s,t\}$} & \multirow{2}{*}{$D_0$} & \multirow{2}{*}{$\alpha$} & \multicolumn{2}{c|}{MILP} & \multicolumn{2}{c|}{CG} & \multicolumn{2}{c}{Greedy} \\
     &  &  & $|C|$ & $t$(s) & $|C|$ & $t$(s) & $|C|$ & $t$(s) \\ \hline
    \multirow{5}{*}{$\{100,9900\}$} & \multirow{5}{*}{4.43} & 15 & 224 & 122 & 224 & 64.7 & 0 & 60.5 \\
      &  & 50 & 134 & 122 & 224 & 63.8 & 0 & 60.5 \\
      &  & 500 & 129 & 123 & 224 & 63.8 & 0 & 60.4 \\
      &  & 1000 & 129 & 122 & 224 & 64 & 0 & 60.4 \\
      &  & 5000 & 138 & 122 & 224 & 64 & 0 & 60.3 \\
    \hline
    \multirow{5}{*}{$\{500,5000\}$} & \multirow{5}{*}{3.30} & 15 & 225 & 122 & 222 & 64.1 & 226 & 8.8 \\
      &  & 50 & 226 & 107 & 222 & 64.3 & 226 & 7.09 \\
      &  & 500 & 223 & 122 & 222 & 62.4 & 226 & 7.01 \\
      &  & 1000 & 223 & 122 & 222 & 64.4 & 226 & 7.13 \\
      &  & 5000 & 223 & 122 & 222 & 63.7 & 226 & 7.1 \\
    \hline
    \multirow{5}{*}{$\{1000,8000\}$} & \multirow{5}{*}{4.06} & 15 & 225 & 122 & 226 & 62.8 & 225 & 21.3 \\
      &  & 50 & 25 & 122 & 226 & 62.4 & 225 & 25.6 \\
      &  & 500 & 22 & 122 & 226 & 63.4 & 225 & 26 \\
      &  & 1000 & 18 & 122 & 226 & 63.1 & 225 & 25.8 \\
      &  & 5000 & 12 & 122 & 226 & 65.3 & 225 & 25.5 \\
    \hline
    \multirow{5}{*}{$\{2500,7500\}$} & \multirow{5}{*}{3.84} & 15 & 179 & 122 & 232 & 50 & 232 & 32.6 \\
      &  & 50 & 186 & 122 & 232 & 49.8 & 232 & 32.5 \\
      &  & 500 & 17 & 122 & 232 & 50.7 & 232 & 32.4 \\
      &  & 1000 & 19 & 122 & 232 & 50.7 & 232 & 32.6 \\
      &  & 5000 & 22 & 123 & 232 & 41 & 232 & 32.8 \\
    \hline
    \end{tabular}
\end{table}

Across the dense instances the ranking is therefore unambiguous: column generation is the only method that is at once fast (about $3$--$7$ seconds with the base set and up to about $67$ seconds with the larger one) and consistently among the best in solution quality, the greedy heuristic is a correct but increasingly expensive fallback that can itself time out, and the exact MILP is no longer usable.

\subsubsection{Sparse 15k Network}
The larger sparse graph at $15{,}000$ nodes confirms the trend; we report only the edge-weighted case (Tables~\ref{tab:syn_sparse15k_w_b120} and~\ref{tab:syn_sparse15k_w_b180}). With the base causality set (Table~\ref{tab:syn_sparse15k_w_b120}) the MILP solves most pairs to optimality, in $17$--$122$ seconds, and the three methods then agree. The exception is $s=100,\,t=9900$ at $\alpha=1.5$, where the MILP runs into its time limit and stops at $150$ causalities (as does the greedy) while column generation reaches $154$; the gap closes at larger $\alpha$, where all methods agree at $155$. Column generation is by far the fastest, often well below one second (e.g.,\ $0.23$ seconds for $s=2500,\,t=7500$), and matches the best count throughout; the greedy trails the best by four for $s=100,\,t=9900$ at $\alpha=1.5$ ($150$ against $154$) and by one for $s=1000,\,t=8000$ at $\alpha=1.5$ ($156$ against $157$), matching it elsewhere.

\begin{table}[!t]
    \centering
    \caption{Synthetic sparse 15k network, edge-weighted (bound $=120$, $h=160$)}
    \label{tab:syn_sparse15k_w_b120}
    \footnotesize
    \setlength{\tabcolsep}{3pt}
    \begin{tabular}{c|c|c|cc|cc|cc}
    \hline
    \multirow{2}{*}{$\{s,t\}$} & \multirow{2}{*}{$D_0$} & \multirow{2}{*}{$\alpha$} & \multicolumn{2}{c|}{MILP} & \multicolumn{2}{c|}{CG} & \multicolumn{2}{c}{Greedy} \\
     &  &  & $|C|$ & $t$(s) & $|C|$ & $t$(s) & $|C|$ & $t$(s) \\ \hline
    \multirow{3}{*}{$\{100,9900\}$} & \multirow{3}{*}{10.34} & 1.5 & 150 & 122 & 154 & 1.96 & 150 & 14 \\
      &  & 2.5 & 155 & 92.5 & 155 & 5.32 & 155 & 13.5 \\
      &  & 3.5 & 155 & 68.7 & 155 & 5.22 & 155 & 14.1 \\
    \hline
    \multirow{3}{*}{$\{500,5000\}$} & \multirow{3}{*}{7.58} & 1.5 & 154 & 119 & 154 & 0.382 & 154 & 3.8 \\
      &  & 2.5 & 155 & 91.6 & 155 & 2.37 & 155 & 3.59 \\
      &  & 3.5 & 155 & 36.7 & 155 & 3.08 & 155 & 2.92 \\
    \hline
    \multirow{3}{*}{$\{1000,8000\}$} & \multirow{3}{*}{6.99} & 1.5 & 157 & 112 & 157 & 0.185 & 156 & 1.05 \\
      &  & 2.5 & 159 & 38.9 & 159 & 1.62 & 159 & 1.07 \\
      &  & 3.5 & 159 & 43.2 & 159 & 2.8 & 159 & 1.05 \\
    \hline
    \multirow{3}{*}{$\{2500,7500\}$} & \multirow{3}{*}{10.49} & 1.5 & 159 & 17.4 & 159 & 0.233 & 159 & 11 \\
      &  & 2.5 & 159 & 17 & 159 & 0.253 & 159 & 11.2 \\
      &  & 3.5 & 159 & 19 & 159 & 0.307 & 159 & 10.9 \\
    \hline
    \end{tabular}
\end{table}

With the harder set (Table~\ref{tab:syn_sparse15k_w_b180}) the MILP collapses again, running to its limit of about $121$ seconds and returning values scattered between $19$ and $227$ causalities. The most dramatic case is $s=1000,\,t=8000$ at $\alpha=1.5$, where it applies only $19$ causalities while the greedy finds $226$ and column generation $228$. Column generation returns the largest causality count on all twelve configurations, between $224$ and $231$, in anything from a fraction of a second up to about a minute; the greedy is a reasonable runner-up, several causalities behind and far more stable than the MILP.

\begin{table}[!t]
    \centering
    \caption{Synthetic sparse 15k network, edge-weighted (bound $=180$, $h=240$)}
    \label{tab:syn_sparse15k_w_b180}
    \footnotesize
    \setlength{\tabcolsep}{3pt}
    \begin{tabular}{c|c|c|cc|cc|cc}
    \hline
    \multirow{2}{*}{$\{s,t\}$} & \multirow{2}{*}{$D_0$} & \multirow{2}{*}{$\alpha$} & \multicolumn{2}{c|}{MILP} & \multicolumn{2}{c|}{CG} & \multicolumn{2}{c}{Greedy} \\
     &  &  & $|C|$ & $t$(s) & $|C|$ & $t$(s) & $|C|$ & $t$(s) \\ \hline
    \multirow{3}{*}{$\{100,9900\}$} & \multirow{3}{*}{10.34} & 1.5 & 170 & 121 & 224 & 25 & 221 & 19.6 \\
      &  & 2.5 & 40 & 121 & 226 & 60.3 & 221 & 20.5 \\
      &  & 3.5 & 216 & 121 & 226 & 60.2 & 221 & 19.7 \\
    \hline
    \multirow{3}{*}{$\{500,5000\}$} & \multirow{3}{*}{7.58} & 1.5 & 198 & 122 & 228 & 0.878 & 227 & 4.59 \\
      &  & 2.5 & 227 & 122 & 231 & 40 & 227 & 5.94 \\
      &  & 3.5 & 95 & 121 & 231 & 35.8 & 227 & 5.89 \\
    \hline
    \multirow{3}{*}{$\{1000,8000\}$} & \multirow{3}{*}{6.99} & 1.5 & 19 & 121 & 228 & 0.289 & 226 & 1.69 \\
      &  & 2.5 & 39 & 121 & 228 & 31.3 & 226 & 1.71 \\
      &  & 3.5 & 63 & 121 & 228 & 21 & 226 & 1.69 \\
    \hline
    \multirow{3}{*}{$\{2500,7500\}$} & \multirow{3}{*}{10.49} & 1.5 & 40 & 122 & 227 & 3.13 & 217 & 17.3 \\
      &  & 2.5 & 46 & 121 & 229 & 56.3 & 217 & 21.8 \\
      &  & 3.5 & 61 & 121 & 229 & 60.4 & 217 & 20.8 \\
    \hline
    \end{tabular}
\end{table}

\subsubsection{Discussion}
Across all synthetic networks, one result is consistent: the exact MILP does not scale. On the dense graph and on the larger causality sets it exhausts its time budget and returns an order of magnitude fewer causalities than column generation, which matches the certified optimum wherever the exact model provides one and is otherwise the only method that is at once fast and consistently high in solution quality, in a fraction of the time. The greedy heuristic remains a useful fast fallback, except on the densest instances, where it becomes expensive or even times out. These results motivate column generation as the method of choice for large networks.

\section{Conclusions and Future Work}
\label{S:7}
We studied the length-constrained shortest-path problem under cascading causal failures introduced by Zhang et al.~\cite{zzy_causal}.
We present three methods.
The first is a mixed-integer linear programming formulation with explicit node-integrity variables that decouple the node state from routing.
The second is a column-generation algorithm whose pricing subproblem is solved by a three-level cascade of beam search, label-capped $A^{*}$, and exact $A^{*}$ with warm-start.
The third is a greedy heuristic as a fast baseline.
On the IEEE~300 and Power~5000 benchmarks, the three methods are competitive.
In dense synthetic networks of $10{,}000$ to $15{,}000$ nodes, the picture changes.
The exact MILP either exhausts its time budget or returns an order of magnitude fewer causalities than column generation.
The greedy heuristic remains correct but becomes expensive and can itself time out on the densest instances.
Column generation is the only method that stays both fast and consistently high in solution quality, in a fraction of the time taken by the MILP.
Column generation is therefore the method of choice for large instances of the problem.
The greedy heuristic is a useful fast alternative when a small loss of optimality is tolerable, and the exact MILP remains valuable for small networks where certified optima are required.

For future work, we will look at the impact of causal node failures on $st$-reliability in a network.
That is, considering the probability that each node is functioning or not, we will construct the model of probabilistic causal failures and investigate the probability that the source node $s$ is connected to the terminal node $t$, that is, $st$-reliability, in the presence of causal failures. We will explore other measures of network vulnerability beyond a topological measure such as shortest path, including flow-based measures. Further, we will explore stochastic causal failures, where a probability governs the causal failure of a node set.

\bibliographystyle{unsrt}
\bibliography{reference}

\end{document}